\documentclass[11pt]{article}

\usepackage{amscd,amsmath, amssymb, fancyhdr, epsfig}
\usepackage[matrix,arrow]{xy}

\numberwithin{equation}{section}

\newcommand{\version}{version 1.0,\ \   12.01.2008}

\def\eqref#1{(\ref{#1})}

\newcommand{\arrow}{{\:\longrightarrow\:}}

\newcommand{\C}{{\Bbb C}}
\newcommand{\R}{{\Bbb R}}

\newcommand{\6}{\partial}
\def\1{\sqrt{-1}\:}
\newcommand{\restrict}[1]{{\left|_{{\phantom{|}\!\!}_{#1}}\right.}}
\newcommand{\cntrct}                
{\hspace{2pt}\raisebox{1pt}{\text{$\lrcorner$}}\hspace{2pt}}

\newcommand{\calo}{{\cal O}}

\renewcommand{\tilde}{\widetilde}
\renewcommand{\bar}{\overline}
\renewcommand{\phi}{\varphi}
\renewcommand{\epsilon}{\varepsilon}
\renewcommand{\geq}{\geqslant}
\renewcommand{\leq}{\leqslant}

\newcommand{\End}{\operatorname{End}}

\newcommand{\Id}{\operatorname{Id}}
\newcommand{\id}{\operatorname{\text{\sf id}}}
\newcommand{\Vol}{\operatorname{Vol}}
\newcommand{\Hom}{\operatorname{Hom}}
\newcommand{\Hol}{\operatorname{Hol}}

\newcommand{\Sing}{\operatorname{Sing}}
\newcommand{\codim}{\operatorname{codim}}

\newcounter{Mycounter}[section]
\newcounter{lemma}[section]
\renewcommand{\thelemma}{{Lemma \thesection.\arabic{lemma}}}
\newcommand{\lemma}{%
    \setcounter{lemma}{\value{Mycounter}}
    \refstepcounter{lemma}
    \stepcounter{Mycounter}
    {\noindent \bf \thelemma:\ }}

\newcounter{claim}[section]
\renewcommand{\theclaim}{{Claim \thesection.\arabic{claim}}}
\newcommand{\claim}{%
    \setcounter{claim}{\value{Mycounter}}
    \refstepcounter{claim}
    \stepcounter{Mycounter}
    {\noindent \bf \theclaim:\ }}

\newcounter{sublemma}[section]
\newcounter{corollary}[section]
\newcounter{theorem}[section]
\renewcommand{\thetheorem}{{Theorem \thesection.\arabic{theorem}}}
\newcommand{\theorem}{%
    \setcounter{theorem}{\value{Mycounter}}
    \refstepcounter{theorem}
    \stepcounter{Mycounter}
    {\noindent \bf \thetheorem:\ }}

\newcounter{conjecture}[section]
\newcounter{proposition}[section]
\renewcommand{\theproposition}
      {{Proposition \thesection.\arabic{proposition}}}
\newcommand{\proposition}{%
    \setcounter{proposition}{\value{Mycounter}}
    \refstepcounter{proposition}
    \stepcounter{Mycounter}
    {\noindent \bf \theproposition:\ }}

\newcounter{definition}[section]
\renewcommand{\thedefinition}
      {{Definition~\thesection.\arabic{definition}}}
\newcommand{\definition}{%
    \setcounter{definition}{\value{Mycounter}}
    \refstepcounter{definition}
    \stepcounter{Mycounter}
    {\noindent \bf \thedefinition:\ }}

\newcounter{example}[section]
\newcounter{remark}[section]
\renewcommand{\theremark}{{Remark \thesection.\arabic{remark}}}
\newcommand{\remark}{%
    \setcounter{remark}{\value{Mycounter}}
    \refstepcounter{remark}
    \stepcounter{Mycounter}
    {\noindent \bf \theremark:\ }}

\newcounter{problem}[section]
\newcounter{question}[section]
\makeatletter

\@addtoreset{equation}{section} \@addtoreset{footnote}{section}
\makeatother

\def\blacksquare{\hbox{\vrule width 5pt height 5pt depth 0pt}}
\def\endproof{\blacksquare}

\begin{document}
\begin{center}
{\LARGE\bf
Positive forms on hyperk\"ahler manifolds\\[4mm]
}

 Misha
Verbitsky\footnote{Misha Verbitsky is 
supported by CRDF grant RM1-2354-MO02.}

\end{center}

{\small \hspace{0.15\linewidth}
\begin{minipage}[t]{0.75\linewidth}
{\bf Abstract} \\ 
Let $(M,I,J,K,g)$ be a hyperk\"ahler manifold,
$\dim_\R M =4n$. We study positive, $\6$-closed 
$(2p,0)$-forms on $(M,I)$.  These forms are 
quaternionic analogues of the positive $(p,p)$-forms, 
well-known in complex geometry. We construct a monomorphism
${\cal V}_{p,p}:\;
  \Lambda^{2p,0}_I(M)\arrow\Lambda^{n+p, n+p}_I(M)$,
which maps $\6$-closed $(2p,0)$-forms to 
closed $(n+p,n+p)$-forms, and positive $(2p,0)$-forms
to positive $(n+p,n+p)$-forms. This construction
is used to prove a hyperk\"ahler version of the 
classical Skoda-El Mir theorem, which says that 
a trivial extension of a closed, positive current
over a pluripolar set is again closed. We also
prove the hyperk\"ahler version of the Sibony's lemma,
showing that a closed, positive $(2p,0)$-form
defined outside of a compact complex subvariety $Z\subset (M,I)$,
$\codim Z > 2p$ is locally integrable in a neighbourhood 
of $Z$. These results are used to prove polystability of
derived direct images of certain coherent sheaves.

\end{minipage}
}

\tableofcontents


\section{Introduction}


\subsection{Hypercomplex manifolds and hyperk\"ahler manifolds}
\label{_hc_hk_Intro_Subsection_}

Let $M$ be a smooth manifold, equipped with an action
of the algebra
\[ {\Bbb H}= \bigg \langle 1, I,J,K\ \ \ | \ \ \ I^2=J^2=IJK=-1\bigg\rangle
\] of quaternions on its tangent
bundle. Such a manifold is called {\bf an almost
  hypercomplex manifold}.
The operators $I$, $J$, $K$ define three almost complex 
structures on $M$. If these almost complex structures
are integrable, $(M,I,J,K)$ is called {\bf a hypercomplex
manifold}.

Hypercomplex manifolds can be defined in terms of
complex geometry, using the notion of a twistor space
(\cite{_HKLR_}, \cite{_Verbitsky:hypercomple_}). 
A scheme-theoretic definition of a hypercomplex
space also exists, allowing one to define hypercomplex
varieties, and even hypercomplex schemes
(\cite{_Verbitsky:hypercomple_}).

Still, in algebraic geometry, the notion of a hyperk\"ahler
manifold is much more popular. A hyperk\"ahler manifold
is a hypercomplex manifold $(M,I,J,K)$, equipped with a Riemannian
form $g$, in such a way that $g$ is a K\"ahler metric with respect
to $I$, $J$ and $K$.

Historically, these definitions were given in opposite
order: Calabi defined the hyperk\"ahler structure in 1978,
and constructed one on the total space of a cotangent bundle
to $\C P^n$ (\cite{_Calabi_}), and Boyer defined
hypercomplex structures and classified compact
hypercomplex manifolds in quaternionic dimension 1
in 1988 (\cite{_Boyer_}). The hyperk\"ahler structures 
are much more prominent because of Calabi-Yau theorem,
\cite{_Yau:Calabi-Yau_},
which can be used to construct hyperk\"ahler structures
on compact, holomorphically symplectic K\"ahler manifolds
(\cite{_Besse:Einst_Manifo_}).

Let $(M,I,J,K,g)$ be a hyperk\"ahler manifold. Since $g$ is 
K\"ahler with respect to $I$, $J$, $K$, the manifold $M$
is equipped with three symplectic forms:
\[
\omega_I(\cdot, \cdot):= g(\cdot, I\cdot), \ \ 
\omega_J(\cdot, \cdot):= g(\cdot, J\cdot), \ \ 
\omega_K(\cdot, \cdot):= g(\cdot, K\cdot).
\]
A simple linear-algebraic calculation can be used to
show that the form $\Omega:= \omega_J + \1\omega_K$
is of Hodge type $(2,0)$ with respect to the complex 
structure $I$ (see e.g. \cite{_Besse:Einst_Manifo_}). 
Since $\Omega$ is also closed, it
is holomorphic. This gives a holomorphic symplectic
structure on a given hyperk\"ahler manifold.
Conversely, each holomorphically symplectic,
compact, K\"ahler manifold admits a hyperk\"ahler
metric, which is unique in a given
K\"ahler class (\cite{_Besse:Einst_Manifo_}).

In algebraic geometry, the words ``hyperk\"ahler''
and ``holomorphically symplectic'' are used as synonyms,
if applied to projective manifolds. There are papers
on ``hyperk\"ahler manifolds in characteristic $p$''
dealing with holomorphically symplectic, projective
manifolds in characteristic $p$.

The first occurence of hyperk\"ahler manifolds precedes
the definition given by Calabi by almost 25 years. 
In his work on classification of irreducible holonomy groups on Riemannian
manifolds, \cite{_Berger:holonomies_}, M. Berger 
listed, among other groups, the group of $Sp(n)$ of
quaternionic unitary matrices. The holonomy of 
the Levi-Civita connection of a K\"ahler manifold
preserves its complex structure (this is one
of the definitions of a K\"ahler manifold). Therefore,
the holonomy of a hyperk\"ahler manifold preserves
$I,J,$ and $K$. We obtained that the holonomy group
of a hyperk\"ahler manifold lies in $Sp(n)$. The converse
is also true: if the Levi-Civita connection of a
Riemannian manifold $M$ preserves a complex structure,
it is K\"ahler (this is, again, one
of the definitions of a K\"ahler manifold),
and if it preserves a action of quaternions,
it is hyperk\"ahler.

In physics, this is often used as a definition
of a hyperk\"ahler structure. 

Summarizing, there are three competing approaches to 
hyperk\"ahler geometry. 
\begin{description}
\item[(i)] A hyperk\"ahler manifold is a Riemannian
manifold $(M,g)$ equipped with almost complex structures $I, J, K$
satisfying $I\circ J = - J\circ I =K$, such that
$(M,I,g)$, $(M,J,g)$ and $(M,K,g)$ are K\"ahler.

\item[(ii)] A hyperk\"ahler manifold is a Riemannian
manifold with holonomy which is a subgroup of $Sp(n)$.

\item[(iii)] (for compact manifolds) 
A hyperk\"ahler manifold is a compact complex manifold 
of K\"ahler type, equipped with a holomorphic symplectic
structure.
\end{description}

Returning to hypercomplex geometry, there is no 
hypercomplex analogue of
Calabi-Yau theorem, hence no definition in terms of
algebro-geometric data such as in (iii).
However, hypercomplex manifold can also be characterized
in terms of holonomy.

Recall that {\em Obata connection} on an almost hypercomplex
manifold is a torsion-free connection preserving $I, J$ and $K$.
Obata (\cite{_Obata_}) has shown that such a connection is
unique, and exists if the almost
complex structures $I$, $J$ and $K$ are integrable.
The holonomy of Obata connection obviously lies
in $GL(n, {\Bbb H})$. The converse is also true:
if a manifold $M$ admits a torsion-free connection
preserving operators $I, J, K\in \End(TM)$,
generating the quaternionic action,
\[ I^2=J^2=K^2 = IJK= - \Id_{TM}, \]
then the almost complex structure operators
$I, J, K$ are integrable. 
Indeed, an almost complex structure is integrable if it is
preserved by some torsion-free connection.

We obtain that a hypercomplex manifold is a manifold
equipped with a torsion-free connection 
$\nabla$ with holonomy $\Hol(\nabla)\subset GL(n, {\Bbb  H})$.
If, in addition, the holonomy of Obata connection is 
a compact group, $M$ is hyperk\"ahler.

Some notions
of complex geometry have  natural quaternionic
analogues in hypercomplex geometry, many of them quite
useful.

By far, the most useful of these is the notion of
HKT-forms, which is a quaternionic analogue of 
K\"ahler forms (\cite{_Gra_Poon_}, \cite{_Banos_Swann_},
\cite{_Alesker_Verbitsky_HKT_}). Generalizing HKT-forms,
one naturally comes across the notion of 
closed, positive $(2,0)$-forms on a hypercomplex
manifold.

\subsection{Positive $(2,0)$-forms on hypercomplex
  manifolds \\ and quaternionic Hermitian structures}

Let $(M,I,J,K)$ be a hypercomplex manifold. We denote the
space of $(p,q)$-forms on $(M,I)$ by $\Lambda^{p,q}_I(M)$.
The operators $I$ and $J$ anticommute, and therefore,
$J(\Lambda^{p,q}_I(M)) = \Lambda^{q,p}_I(M)$. The map
$\eta \arrow J(\bar\eta)$ induces an anticomplex
endomorphism of $\Lambda^{p,q}_I(M)$. Clearly,
\[
{J^2}\restrict{\Lambda^{p,q}_I(M)} = (-1)^{p+q}\Id.
\]
For $p+q$ even, $\eta \arrow J(\bar\eta)$
is an anticomplex involution, that is, a real structure
on $\Lambda^{q,p}_I(M)$. A $(2p,0)$-form 
$\eta\in\Lambda^{2p,0}_I(M)$ is called {\bf real} if 
$\eta=J(\bar\eta)$. The bundle of real $(2p,0)$-forms
is denoted $\Lambda^{2p,0}_I(M, \R)$.

The real $(2,0)$-forms are most significant,
because they can be interpreted as quaternionic 
pseudo-Hermitian structures.

Recall that a Riemannian metric $g$ on an almost complex
manifold $(M,I)$ is called {\bf Hermitian} if $g$ is $U(1)$-invariant,
with respect to the $U(1)$-action on $TM$ defined by 
\[
t \arrow \cos t \cdot \id_{TM} + \sin t \cdot I.
\]
This is equivalent to $g(I\cdot, I\cdot)= g(\cdot, \cdot).$

When $M$ is almost hypercomplex, it is natural to consider
a group $G\subset \End(TM)$ generated by $U(1)$-action associated with
$I$, $J$, $K$ as above. It is easy to see that $G$ is the group
of unitary quaternions, naturally identified with $SU(2)$. 
Thus obtained action of $SU(2)$ on $\Lambda^*(M)$ is fundamental,
and plays in hypercomplex and hyperk\"ahler geometry the same
role as played by the Hodge structures in complex algebraic 
geometry. 

Recall that bilinear symmetric forms (not necessarily
positive definite) on $TM$ are called 
{\em pseudo-Riemannian structures.}

A (pseudo-)Riemannian structure $g$ on an almost
hypercomplex manifold $(M,I,J,K)$ is called
{\bf quaternionic (pseudo-)Hermitian} if $g$
is $SU(2)$-invariant. In other words, a 
quaternionic pseudo-Hermitian structure is 
a bilinear, symmetric, $SU(2)$-invariant form on $TM$.

Given a real $(2,0)$-form $\eta\in \Lambda^{2,0}_I(M,\R)$,
consider a bilinear form \[ g_\eta(x,y):= \eta(x, Jy)\] on $TM$.
Since $\eta$ is a $(2,0)$-form, we have \[ \eta(Ix, Iy) =
- \eta(x,y),\] for all $x, y \in TM$ 
and therefore \[ g_\eta(Ix, Iy)= g_\eta(x,y).\] Similarly, 
we obtain $g_\eta(Jx, Jy)= g_\eta(x,y),$ because 
$\eta\left (J(\bar x), J(\bar y)\right)= \eta(x,y).$

 Since $\eta$ is skew-symmetric, and $J^2=-1$,  
$g_\eta$ is symmetric. We obtained that $g_\eta$
is a pseudo-Hermitian form on $TM$. This construction
is invertible (see Section \ref{_quat_pseudo-H_Section_}),
and gives an isomorphism between the bundle  $H$ of
real $(2,0)$-forms and the bundle $\Lambda^{2,0}_I(M,\R)$ of
quaternionic pseudo-Hermitian forms (\ref{_pseudo-Hermi_2,0_Claim_}).
The inverse isomorphism  $H \arrow \Lambda^{2,0}_I(M,\R)$
is given as follows. Starting from a 
quaternionic pseudo-Hermitian form $g$, we construct
2-forms $\omega_I, \omega_J, \omega_K$ as in
Subsection \ref{_hc_hk_Intro_Subsection_}.
Then $\Omega_g:= \omega_J+ \1\omega_K$ 
is a real $(2,0)$-form. 

A real $(2,0)$-form $\eta$ is called {\bf positive
  definite} if the corresponding symmetric form $g_\eta$
is positive definite. 

There are two differentials on $\Lambda^{*,0}_I(M)$: 
the standard Dolbeault differential
$\6:\; \Lambda^{p,0}_I(M)\arrow \Lambda^{p+1,0}_I(M),$ and $\6_J$, which
is obtained from $\6$ by twisting with $\eta \arrow
J(\bar\eta)$. One could define $\6_J$ as 
$\6_J(\eta):= -J\bar\6(J\eta)$.

The pair of differentials $\6, \6_J$ behaves in many ways
similarly to the operators $d, d^c$ on a complex
manifold. They anticommute, and satisfy $\6^2=\6_J^2=0$.

A positive definite  $(2,0)$-form $\eta \in\Lambda^{2,0}_I(M,\R)$
is called {\bf HKT-form} if $\6\eta=0$. The corresponding
quaternionic Hermitian metric is called {\bf the HKT-metric}.
This notion was first defined by string physicists
\cite{_Howe_Papado_}, and much studied since then (see 
\cite{_Gra_Poon_} for an excellent survey
of an early research). 

In \cite{_Banos_Swann_} 
(see also \cite{_Alesker_Verbitsky_HKT_}),
it was shown that HKT-forms locally always have
a real-valued potential $\phi$, known as HKT-potential:
$\eta = \6\6_J \phi$. This function is a
quaternionic analogue of the K\"ahler potential.

We obtain the following dictionary of
parallels between the complex and hypercomplex manifolds.

\hfill

\begin{tabular}{|c|c|}
\hline
\it $\C$  & \it ${\Bbb H}$\\[1mm]
\hline
complex manifold & hypercomplex manifold \\[1mm]
\hline
$\Lambda^{p,p}(M, \R)$ & $\Lambda^{2p,0}_I(M, \R)$\\[1mm]
\hline
$d, d^c$ & $\6, \6_J$ \\[1mm]
\hline
real $(1,1)$-forms & real $(2,0)$-forms \\[1mm]
\hline
closed positive definite
(1,1)-forms & HKT-forms \\[1mm] 
\hline
K\"ahler potentials & HKT-potentials \\[1mm] \hline
\end{tabular}

\hfill

\noindent
This analogy can be built upon, to obtain the
notion of positive $(2p,0)$-forms.

\subsection{Positive $(2p,0)$-forms on hypercomplex manifolds}

\definition
(\cite{_Alesker_Verbitsky_HKT_})
A real $(2p,0)$-form $\eta\in \Lambda^{2,0}_I(M,\R)$ on a
hypercomplex manifold is called {\bf weakly positive} if
\[
\eta
\left(x_1, J(\bar x_1), x_2, J(\bar x_2), ... x_p, J(\bar
x_p) \right) \geq 0,
\]
for any $x_1, ..., x_p \in T^{1,0}_IM$,
and {\bf closed} if $\6\eta=0$.

\hfill

In modern complex geometry, the positive, closed $(p,p)$
forms and currents play a central role, due to several 
by now classical theorems, which were proven in 1960-1980-ies,
building upon the ideas of P. Lelong (see \cite{_Demailly:L^2_} for
an elementary exposition of the theory of positive currents). 

The hypercomplex analogue of these results could be just
as significant. 

\hfill

In \cite{_Alesker_Verbitsky_HKT_}, 
a hypercomplex version of the classical Chern-Levine-Nirenberg
theorem was obtained. In the present paper, we prove quaternionic
versions of two classical theorems, both of them quite important in
complex geometry.

\hfill

\theorem\label{_Sibony_hk_Theorem_}
(``Sibony's Lemma'')
Let $(M,I,J,K, g)$ be a hyperk\"ahler manifold,
$\dim_\R M=4n$, and $Z\subset (M,I)$ a compact
complex subvariety, $\codim Z > 2p$. Consider
a weakly positive, closed form $\eta \in \Lambda^{2p,0}_I(M\backslash Z, \R)$.
Then $\eta$ is locally integrable around $Z$.

\hfill

{\bf Proof:} See 
\ref{_L^1_inte_Theorem_}. \endproof

\hfill

The classical version of this theorem states that
a closed, positive $(p,p)$-form defined outside
of a complex subvariety of codimension $>p$ is integrable
in a neighbourhood of this subvariety. Its proof can be
obtained by slicing.

In hypercomplex geometry, the slicing is possible only
on a flat manifold, because a typical
hypercomplex manifold has no non-trivial hypercomplex
subvarieties, even locally. In earlier versions of \cite{_V:reflexive_},
 \ref{_Sibony_hk_Theorem_} was proven for flat
 hypercomplex manifold using slicing, and then
extended to non-flat manifold by approximation.
The approximation argument was very unclear
and ugly. In 2007, a new proof of Sibony's lemma
was found (\cite{_Verbitsky:omega-psh_}),
using the emerging theory of plurisubharmonic 
functions on calibrated manifolds (\cite{_Harvey_Lawson:Psh_},
\cite{_Harvey_Lawson:Dua_}) instead of slicing. 
In Section \ref{_Sibony_Section_}, we adapt this
argument to hyperk\"ahler geometry, obtaining a relatively
simple and clean proof of \ref{_Sibony_hk_Theorem_}.

\hfill

\ref{_Sibony_hk_Theorem_} was used in 
 \cite{_V:reflexive_} to prove results about stability of
certain coherent sheaves on hyperk\"ahler manifolds
(Subsection \ref{_hh_intro_Subsection_}). 
\ref{_Sibony_hk_Theorem_} was used to show that
the form representing $c_1(F)$ for such a sheaf 
is integrable. To prove theorems about stability,
we need also to show that the corresponding current
is closed. Then the integral of the form
representing $c_1(F)$ can be interpreted in terms
of the cohomology. 

Given a form $\eta$ on $M \backslash Z$, locally
integrable everywhere on $M$, we can interpret
$\eta$ as a current on $M$,
\[ \alpha \arrow \int_{M \backslash Z} \eta \wedge \alpha.
\]
This current is called {\bf a trivial extension of
$\eta$ to $M$}. A priori, it can be non-closed.
However, in complex geometry, a trivial extension
of an integrable, closed and positive form is again
closed. This fundamental result is known as 
Skoda-El Mir theorem (\ref{_Sk_E-M_complex_geo_Theorem_}).
In Section \ref{_S_E-M_Section_}, we prove a hypercomplex
analogue of Skoda-El Mir theorem.

Recall that hypercomplex manifolds can be defined
in terms of holonomy (Subsection \ref{_hc_hk_Intro_Subsection_}),
as manifolds equipped with a torsion-free connection $\nabla$,
with $\Hol(\nabla) \subset GL(n, {\Bbb H})$.
A hypercomplex manifold $(M,I,J,K)$ is called
{\bf an $SL(n, {\Bbb H})$-manifold} if its
holonomy lies in $SL(n, {\Bbb H})\subset GL(n, {\Bbb H})$.
Such manifolds were studied in \cite{_Verbitsky:canoni_} 
and \cite{_BDV:nilmanifolds_}. It was shown that
$(M,I,J,K)$ is an $SL(n, {\Bbb H})$-manifold if and
only if $M$ admits a holomorphic, real $(2n,0)$-form.
In particular, all hyperk\"ahler manifolds 
satisfy $\Hol(\nabla) \subset SL(n, {\Bbb H})$.

\hfill

\theorem\label{_hk_S-EM_intro_Theorem_}
Let $(M,I,J,K)$ be an $SL(n, {\Bbb H})$-manifold,
and $Z\subset (M,I)$ a closed complex subvariety.
Consider a closed, positive form
\[ \eta \in \Lambda^{2p,0}_I(M\backslash Z, \R),\] and assume that
$\eta$ is locally integrable around $Z$. Let
$\tilde \eta$ be the current obtained as a 
trivial extension of $\eta$ to $M$. Then
$\6\tilde\eta=0$.

\hfill

{\bf Proof:} \ref{_hk_S-EM_Theorem_}. \endproof

\hfill

\ref{_hk_S-EM_intro_Theorem_}
is deduced from the classical Skoda-El Mir theorem.
In Subsection \ref{_V_p,q_Subsection_},
we construct a map ${\cal V}_{p,q}:\;
  \Lambda^{p+q,0}_I(M)\arrow\Lambda^{n+p, n+q}_I(M)$,
which has the following properties.

\hfill

\claim\label{_V_properties_Claim_}
Let $\eta\in \Lambda^{2p,0}_I(M)$ be a $(2p,0)$-form
on an $SL(n, {\Bbb H})$-manifold.
Then the $(n+p, n+p)$-form $(\1)^p {\cal V}_{p,p}(\eta)$ is 
real (in the usual sense) if and only if $\eta$ is real,
positive if and only if $\eta$ is positive, and closed
if and only if $\6_J\eta = \6\eta=0$. 

\hfill

{\bf Proof:} Follows immediately from 
\ref{_V_main_Proposition_}. \endproof

\hfill

To prove \ref{_hk_S-EM_intro_Theorem_},
take $\eta \in \Lambda^{2p,0}_I(M)$ which is closed
and positive. As follows from \ref{_V_properties_Claim_},
the $(n+p,n+p)$-form $(\1)^p {\cal V}_{p,p}(\eta)$ 
is closed and positive, in the usual complex-analytic
sense. Its trivial extension is closed and positive, 
by the Skoda-El Mir theorem. Then $(\1)^p {\cal V}_{p,p}(\tilde \eta)$ 
is closed. Applying \ref{_V_properties_Claim_}
again, we find that closedness of 
$(\1)^p {\cal V}_{p,p}(\tilde \eta)$ 
implies that $\6\tilde \eta=0$.

\subsection{Hyperholomorphic bundles and reflexive sheaves}
\label{_hh_intro_Subsection_}

The results about positive $(2,0)$-forms on hypercomplex
manifolds are especially useful in hyperk\"ahler geometry.
In \cite{_V:reflexive_}, we used this notion to prove theorems
about stability of direct images of coherent sheaves. The
earlier arguments were unclear and flawed, and the machinery of positive
$(2p,0)$-forms was developed in order to obtain clear proofs
of these results. Here we give a short sketch of main arguments
used in \cite{_V:reflexive_}. Throughout this paper, stability
of coherent sheaves is understood in Mumford-Takemoto sense.

Let $(M, I,J,K)$ be a compact hyperk\"ahler manifold, and $B$ a holomorphic
Hermitian bundle on $(M,I)$. Denote the Chern connection on $B$
by $\nabla$. We say that $B$ is {\bf hyperholomorphic}
if its curvature $\Theta_B \in \Lambda^2(M) \otimes \End B$
is $SU(2)$-invariant, with respect to the natural action
of $SU(2)$ on $\Lambda^2(M)$. This notion was defined in 
\cite{_Verbitsky:Hyperholo_bundles_},
and much studied since then.

It is easy to check that $SU(2)$-invariant 2-forms are 
pointwise orthogonal to the K\"ahler form $\omega_I$.
Therefore, $(B, \nabla)$ satisfies the Yang-Mills
equation $\Lambda\Theta_B=0$.\footnote{Here
\[ \Lambda^{p,q}_I(M) \otimes \End B \arrow \Lambda^{p,q}_I(M) \otimes \End B\]
is the standard Hodge operator, which is 
Hermitian adjoint to $L(\eta) = \omega_I \wedge \eta$.}  
In other words, $\nabla$ is Hermitian-Einstein.

One can easily prove that Yang-Mills bundles are always
{\bf polystable}, that is, obtained as a direct sum of stable
bundles of the same slope. The converse is also true:
as follows from Donaldson-Uhlenbeck-Yau theorem \cite{_Uhle_Yau_}, 
a Yang-Mills connection exists on any polystable bundle,
and is unique. 

In \cite{_Verbitsky:Hyperholo_bundles_},
it was shown that a polystable bundle on $(M,I)$
admits a hyperholomorphic connection if and only if the 
Chern classes $c_1(B)$ and $c_2(B)$ are $SU(2)$-invariant.

In \cite{_Verbitsky:Symplectic_II_}, it was shown that
for any compact hyperk\"ahler manifold \\ $(M,I,J,K)$
there exists a countable set 
\[ P\subset S^2 = \{ a, b, c \ \ | \ \ a^2 +b^2 + c^2=1\}
\] 
with the following property. For any
$(a, b, c)\notin P$, let $L:=aI + bJ + cK$ be the
corresponding complex structure on $M$ induced by the 
quaternionic action. Then all integer $(p,p)$-classes
on $(M,L)$ are $SU(2)$-invariant. In particular, all
stable bundles on $(M, L)$ are hyperholomorphic.

Many of these results can be extended to reflexive coherent sheaves.
Recall that a coferent sheaf $F$ on a complex manifold $X$
is called {\bf reflexive} if the natural map $F \arrow F^{**}$
is an isomorphism. Here, $F^*$ denotes the dual sheaf,
$F^*:= \Hom(F, \calo_X)$. The following properties 
of reflexive sheaves are worth mention (see \cite{_OSS_}).

\begin{itemize}
\item Holomorphic vector bundles are obviously reflexive.

\item Let $Z\subset X$ be a closed complex subvariety,
$\codim Z \geq 2$, and $j:\; X \backslash Z \arrow X$
the natural embedding. Then $j_*F$ is reflexive, for
any reflexive sheaf $F$ on $X \backslash Z$.

\item The sheaf $F^*$ is reflexive, for any coherent sheaf
  $F$.

\item For any torsion-free coherent sheaf $F$, the natural map
$F \arrow F^{**}$ is a monomorphism, and $F^{**}$ is reflexive.
Moreover, $F^{**}$ is a minimal reflexive sheaf containing $F$.

\item For any torsion-free coherent sheaf $F$, the 
singular set $\Sing(F)$ has codimension $\geq 2$.
If $F$ is reflexive,  $\Sing(F)$ has codimension $\geq 3$.

\item A torsion-free sheaf of rank 1 is always reflexive.

\item A torsion-free sheaf $F$ is stable if and only if $F^{**}$
is stable.
\end{itemize}

In \cite{_V:Hyperholo_sheaves_}, the definition of a hyperholomorphic
connection was extended to reflexive coherent sheaves, using the notion of
admissible connection defined by Bando and Siu 
in a fundamental work \cite{_Bando_Siu_}.

\hfill

Let us recall what Bando and Siu did.

\hfill

\definition
Let $(X, \omega)$ be a K\"ahler manifold, $Z \subset X$ a
closed complex subvariety, $\codim Z\geq 2$, and $F$ 
a holomorphic vector bundle on $X \backslash Z$. 
Given a Hermitian metric $h$ on $F$, denote by $\nabla$
the corresponding Chern connection, and let $\Theta_F$ be
its curvature. The metric $h$ and the connection $\nabla$ are called
{\bf admissible} if 
\begin{description}
\item[(i)] $\Lambda\Theta_F$ is uniformly bounded, where 
$\Lambda:\; \Lambda^{1,1}_I(M) \otimes \End B\arrow \End B$
is the Hodge operator, which is 
Hermitian adjoint to $L(\eta) = \omega_I \wedge \eta$.
\item[(ii)] The curvature $\Theta_F$ is locally
  $L^2$-integrable everywhere on $M$. 
\end{description}

\hfill

Bando and Siu proved the following.

\begin{itemize}

\item Let $(X, \omega)$ be a K\"ahler manifold, $Z \subset X$ a
closed complex subvariety, $\codim Z\geq 2$, and $F$ 
a holomorphic vector bundle on $X \backslash Z\stackrel \hookrightarrow X$.
Assume that $F$ is equipped with an admissible connection.
Then $j_* F$ is a reflexive coherent sheaf. Conversely,
any coherent sheaf admits an admissible connection
outside of its singularities. Such a connection is called
{\bf an admissible connection on $F$}.

\item A version of Donaldson-Uhlenbeck-Yau theorem
is valid for stable reflexive sheaves. Let $F$ be a reflexive
sheaf on a compact K\"ahler manifold $X$. The 
admissible connection on $F$ is called {\bf Yang-Mills}
if $\Lambda \Theta_F= c \Id_F$, where $\Theta_F$
is its curvature, and $c$ some constant.
Bando and Siu proved that a Yang-Mills
connection is unique, and exists if and only
if $F$ is polystable. 
\end{itemize}

In \cite{_V:Hyperholo_sheaves_}, these results were applied
to coherent sheaves on a hyperk\"ahler manifold $(M,I,J,K, g)$.
A {\bf hyperholomorphic connection} on a reflexive sheaf $F$ on $(M,I)$ 
is an admissible connection with $SU(2)$-invariant curvature.
Since any $SU(2)$-invariant form $\Theta_F$ satisfies 
$\Lambda \Theta_F= 0$, a hyperholomorphic connection 
is always Yang-Mills. In \cite{_V:Hyperholo_sheaves_},
it was shown that any polystable reflexive sheaf with
$SU(2)$-invariant  Chern classes $c_1(F)$, $c_2(F)$ 
admits a hyperholomorphic connection.
 
In \cite{_V:reflexive_}, this formalism 
was used to prove polystability of derived direct images of
hyperholomorphic bundles. Let $M_1, M_2$ be compact
hyperk\"ahler manifolds, and $B$ a hyperholomorphic 
bundle on $M_1\times M_2$. Denote the natural projection
$M_1 \times M_2 \arrow M_2$ by $\pi$. It was shown that
the derived direct image sheaves $R^i\pi_* B$ admit
a hyperholomorphic connection, outside of their
singularities. Were this connection admissible, Bando-Siu
theorem would imply polystability of $R^i\pi_* B$ outright.
However, $L^2$-integrability of its curvature 
is difficult to establish. In  \cite{_V:reflexive_},
we proposed a roundabout argument to prove polystability
of $F:= (R^i\pi_* B)^{**}$.

Let $(M,I,J,K,g)$ be a compact hyperk\"ahler manifold,
$\dim_\R M = 4n$, and $F$ a reflexive coherent sheaf on $(M,I)$. Assume
that outside of its singularities, $F$ is equipped with
a metric, and its Chern connection has $SU(2)$-invariant
curvature. Consider a subsheaf $F_1\subset F$. Then,
outside of singularities of $F$, $F_1$, the class
$-c_1(F)$ is represented by a form $\nu$ with
$\nu - J(\nu)$ positive, and vanishing only
if $F= F_1 \oplus F_2$. This follows from an argument
which is similar to one that proves that holomorphic
subbundles of a flat bundle have negative $c_1$: the
$SU(2)$-invariance of the curvature $\Theta_F$
is equivalent to $\Theta_F - J(\Theta_F)=0$. 
From positivity and non-vanishing of $\nu - J(\nu)$,
one needs to infer that $\deg c_1(F_1) <0$, which 
would suffice to show that $F$ is polystable. 

The expression
\begin{equation}\label{_deg_via_nu_Equation_}
\deg c_1(F_1) = - \int_M \nu\wedge \omega_I^{2n-1}
= -\frac  1 2  \int_M (\nu - J(\nu))\wedge \omega_I^{2n-1} 
\end{equation}
would have been true were the form $\nu- J(\nu)$ integrable, 
and closed as a current on $M$. However, the
$(2,0)$-form $\Omega_\nu$ corresponding to $\nu$ as in
Section \ref{_quat_pseudo-H_Section_}
is $\6$-closed,  because $\nu$ is closed. This form is
positive, because $\nu - J(\nu)$ is positive,
and $\Omega_\nu$ satisfies $2\Omega_\nu = \Omega_{\nu -J\nu}$,
which is clear from its construction. This form
is defined outside of the set $S\subset M$ where the 
sheaves $F, F_1$ are not locally trivial.
Since these sheaves are reflexive, $\codim S>2$,
and we could apply the hyperk\"ahler version of
Sibony's lemma (\ref{_Sibony_hk_Theorem_})
to obtain that $\Omega_\nu$ is  integrable. 
Now, the hypercomplex version of
Skoda-El Mir theorem (\ref{_hk_S-EM_intro_Theorem_})
implies that the trivial extension of
$\Omega_\nu$ is a $\6$-closed current.
Therefore, $\deg F_1$ can be computed 
through the integral \eqref{_deg_via_nu_Equation_}.
Since $\nu - J(\nu)$ is positive, this integral
is negative, and strictly negative unless
$F= F_1 \oplus F_2$. Therefore, $F$ is polystable.
We gave a sketch of an argument showing that
$F= (R^i\pi_* B)^{**}$ is polystable. For a complete
proof, please see \cite{_V:reflexive_}.


\section{Quaternionic Dolbeault complex}


In this Section, we introduce the quaternionic Dolbeault 
complex \[ \left(\bigoplus \Lambda^{p,q}_{I,+}, d_+\right),\] 
used further on in this paper.
We follow \cite{_Verbitsky:HKT_}.

\subsection{Weights of $SU(2)$-representations}
\label{_Weights_Subsection_}

It is well-known that any irreducible representation
of $SU(2)$ over $\C$ can be obtained as a symmetric power
$S^i(V_1)$, where $V_1$ is a fundamental 2-dimensional
representation. We say that a representation $W$ 
{\bf has weight $i$} if it is isomorphic to $S^i(V_1)$.
A representation is said to be {\bf pure of weight $i$}
if all its irreducible components have weight $i$.
If all irreducible components of a representation $W_1$
have weight $\leq i$, we say that $W_1$ {\bf is a
  representation of weight $\leq i$}.
In a similar fashion one defines representations
of weight $\geq i$.

\hfill

\remark\label{_weight_multi_Remark_}
The Clebsch-Gordan formula (see \cite{_Humphreys_})
claims that the weight is {\em multiplicative}, 
in the following sense: if $i\leq j$, then
\[
V_i\otimes V_j = \bigoplus_{k=0}^i V_{i+j-2k},
\]
where $V_i=S^i(V_1)$ denotes the irreducible
representation of weight $i$.

\hfill

A subspace $W\subset W_1$ is {\bf pure of weight $i$}
if the $SU(2)$-representation $W'\subset W_1$ generated
by $W$ is pure of weight $i$.

\subsection{Quaternionic Dolbeault complex: a definition}
\label{_qD_Subsection_}

Let $M$ be a hypercomplex (e.g. a hyperk\"ahler) manifold,
$\dim_{\Bbb H}M=n$.
There is a natural multiplicative action of $SU(2)\subset
{\Bbb H}^*$ on $\Lambda^*(M)$, associated with the
hypercomplex structure. 

\hfill

\remark\label{_weights_on_forms_Remark_}
The space $\Lambda^*(M)$
is an infinite-dimensional representation of $SU(2)$,
however, all its irreducible components
are finite-dimensional. Therefore it makes
sense to speak of {\em weight} of $\Lambda^*(M)$
and its sub-\-rep\-re\-sen\-ta\-tions. Clearly, $\Lambda^1(M)$
has weight 1. From Clebsch-Gordan formula
(\ref{_weight_multi_Remark_}), it follows that 
 $\Lambda^i(M)$ is an $SU(2)$-representation
of weight $\leq i$. Using the Hodge $*$-isomorphism 
$\Lambda^i(M)\cong \Lambda^{4n-i}(M)$, we find that
for $i> 2n$, $\Lambda^i(M)$ is a representation
of weight $\leq 2n-i$.

\hfill

Let $V^i\subset \Lambda^i(M)$ be a maximal
$SU(2)$-invariant subspace of weight $<i$.
The space $V^i$ is well defined, because
it is a sum of all irreducible representations
$W\subset \Lambda^i(M)$ of weight $<i$.
Since the weight is multiplicative
(\ref{_weight_multi_Remark_}), $V^*= \bigoplus_i V^i$
is an ideal in $\Lambda^*(M)$. We also have
$V^i = \Lambda^i(M)$ for $i> 2n$
(\ref{_weights_on_forms_Remark_}).

It is easy to see that the de Rham differential
$d$ increases the weight by 1 at most. Therefore,
$dV^i\subset V^{i+1}$, and $V^*\subset \Lambda^*(M)$
is a differential ideal in the de Rham DG-algebra
$(\Lambda^*(M), d)$.

\hfill

\definition\label{_qD_Definition_}
Denote by $(\Lambda^*_+(M), d_+)$ the quotient algebra
$\Lambda^*(M)/V^*$
It is called {\bf the quaternionic Dolbeault algebra of
  $M$}, or {\bf the quaternionic Dolbeault complex} 
(qD-algebra or qD-complex for short).

The space $\Lambda^i_+(M)$ can be identified with the
maximal subspace of $\Lambda^i(M)$ of weight $i$,
that is, a sum of all irreducible sub-representations of weight $i$.
This way, $\Lambda^i_+(M)$ can be considered as a subspace
in $\Lambda^i(M)$; however, this subspace is not preserved
by the multiplicative structure and the differential.

\hfill

\remark 
The complex $(\Lambda^*_+(M), d_+)$ 
was constructed much earlier by Salamon,
in a different (and much more general) situation,
and much studied since then
(\cite{_Salamon_}, \cite{_Capria-Salamon_},
\cite{_Baston_}, \cite{_Leung_}).

\subsection{The Hodge decomposition of the quaternionic
  Dolbeault complex}.
\label{_Hodge_on_qD_Subsection_}

Let $(M,I,J,K)$ be a hypercomplex 
manifold, and $L$ a complex structure
 induced by the quaternionic action, say,
$I$, $J$ or $K$.
 Consider the $U(1)$-action
on $\Lambda^1(M)$ provided by 
$\phi\stackrel{\rho_L} \arrow \cos \phi\Id + \sin\phi \cdot L$.
We extend this action to a multiplicative action on
$\Lambda^*(M)$. Clearly, for a $(p,q)$-form 
$\eta\in \Lambda^{p,q}(M,L)$, we have
\begin{equation}\label{_Hodge_weights_Equation_}
   \rho_L(\phi)\eta = e^{\1(p-q)\phi}\eta.
\end{equation}

This action is compatible with the weight decomposition
of $\Lambda^*(M)$, and gives a Hodge decomposition 
of $\Lambda^*_+(M)$ (\cite{_Verbitsky:HKT_}).
\[
\Lambda^i_+(M) = \bigoplus_{p+q=i}\Lambda^{p,q}_{+,I}(M)
\]

\hfill

The following result is implied immediately by the
standard calculations from the theory of
$SU(2)$-representations.

\hfill

\proposition \label{_qD_decompo_expli_Proposition_}
Let $(M,I,J,K)$ be a hypercomplex manifold and
\[
\Lambda^i_+(M) = \bigoplus_{p+q=i}\Lambda^{p,q}_{+,I}(M)
\]
the Hodge decomposition of qD-complex defined above.
Then there is a natural isomorphism
\begin{equation}\label{_qD_decompo_Equation_}
\Lambda^{p,q}_{+,I}(M)\cong \Lambda^{p+q,0}(M,I).
\end{equation}

{\bf Proof:} See \cite{_Verbitsky:HKT_}. \endproof

\hfill

This isomorphism is compatible with a natural algebraic
structure on $\bigoplus_{p+q=i}\Lambda^{p+q,0}(M,I)$,
and with the Dolbeault differentials, in the following
way.

\hfill

Let $(M,I,J,K)$ be a hypercomplex manifold.
We extend \[ J:\; \Lambda^1(M) \arrow \Lambda^1(M)\]
to $\Lambda^*(M)$ by multiplicativity. Recall that 
\[ J(\Lambda^{p,q}(M,I))=\Lambda^{q,p}(M,I), \]
because $I$ and $J$ anticommute on $\Lambda^1(M)$.
Denote by 
\[ \6_J:\;  \Lambda^{p,q}(M,I)\arrow \Lambda^{p+1,q}(M,I)
\]
the operator $J\circ \bar\6 \circ J$, where
$\bar\6:\;  \Lambda^{p,q}(M,I)\arrow \Lambda^{p,q+1}(M,I)$
is the standard Dolbeault operator on $(M,I)$, that is, the
$(0.1)$-part of the de Rham differential.
Since $\bar\6^2=0$, we have $\6_J^2=0$.
In \cite{_Verbitsky:HKT_} it was shown that $\6$ and $\6_J$
anticommute:
\begin{equation}\label{_commute_6_J_6_Equation_}
\{\6_J, \6 \}=0.
\end{equation}

Consider the quaternionic Dolbeault complex
$(\Lambda^*_+(M), d_+)$ constructed in Subsection
\ref{_qD_Subsection_}. Using the Hodge decomposition, we can
represent this complex as

\begin{equation}\label{_bicomple_d_+_Equation}
\begin{minipage}[m]{0.85\linewidth}
\begin{center}{ $
\xymatrix @C+1mm @R+10mm@!0  { 
  && \Lambda^0_{+,I}(M) \ar[dl]^{d^{1,0}_{+,I}} \ar[dr]^{d^{0,1}_{+,I}} 
   &&  \\
 & \Lambda^{1,0}_{+,I}(M) \ar[dl]^{d^{1,0}_{+,I}} \ar[dr]^{d^{0,1}_{+,I}} &
 & \Lambda^{0,1}_{+,I}(M) \ar[dl]^{d^{1,0}_{+,I}} \ar[dr]^{d^{0,1}_{+,I}}&\\
 \Lambda^{2,0}_{+,I}(M) && \Lambda^{1,1}_{+,I}(M) 
   && \Lambda^{0,2}_{+,I}(M) \\
}
$
}\end{center}
\end{minipage}
\end{equation}
where $d^{1,0}_{+,I}$,  $d^{0,1}_{+,I}$ are the Hodge components of 
the quaternionic Dolbeault differential $d_+$, taken with
respect to $I$. 

\hfill

\theorem\label{_bico_ide_Theorem_}
Under the isomorphism 
\[
\Lambda^{p,q}_{+,I}(M)\cong \Lambda^{p+q,0}(M,I)
\]
constructed in \ref{_qD_decompo_expli_Proposition_},
$d^{1,0}_+$
corresponds to $\6$ and $d^{0,1}_+$
to $\6_J$:

\begin{equation}\label{_bicomple_XY_Equation}
\begin{minipage}[m]{0.85\linewidth}
{\tiny $
\xymatrix @C+1mm @R+10mm@!0  { 
  && \Lambda^0_+(M) \ar[dl]^{d^{0,1}_+} \ar[dr]^{d^{1,0}_+} 
   && && && \Lambda^{0,0}_I(M) \ar[dl]^{\6} \ar[dr]^{ \6_J}
   &&  \\
 & \Lambda^{1,0}_+(M) \ar[dl]^{d^{0,1}_+} \ar[dr]^{d^{1,0}_+} &
 & \Lambda^{0,1}_+(M) \ar[dl]^{d^{0,1}_+} \ar[dr]^{d^{1,0}_+}&& 
\text{\large $\cong$} &
 &\Lambda^{1,0}_I(M)\ar[dl]^{ \6} \ar[dr]^{ \6_J}&  &
 \Lambda^{1,0}_I(M)\ar[dl]^{ \6} \ar[dr]^{ \6_J}&\\
 \Lambda^{2,0}_+(M) && \Lambda^{1,1}_+(M) 
   && \Lambda^{0,2}_+(M)& \ \ \ \ \ \ & \Lambda^{2,0}_I(M)& & 
\Lambda^{2,0}_I(M) & &\Lambda^{2,0}_I(M) \\
}
$
}
\end{minipage}
\end{equation}

\hfill

{\bf Proof:} See \cite{_Verbitsky:HKT_} or \cite{_Verbitsky:qD_}.
For another proof \ref{_bico_ide_Theorem_},
please see \ref{_R_p,q_and_6_Claim_}. \endproof


\section{Quaternionic pseudo-Hermitian structures}
\label{_quat_pseudo-H_Section_}


Further on in this paper, we shall use 
some results about diagonalization of
certain $(2,0)$-forms associated to 
quaternionic pseudo-Hermitian structures.
The results of this section are purely
linear-algebraic and elementary. We follow 
\cite{_Verbitsky:HKT_}, \cite{_Verbitsky_HKT-exa_} 
and \cite{_Alesker_Verbitsky_HKT_}.

Let $(M,I,J,K)$ be a hypercomplex manifold.
A quaternionic pseudo-Hermitian form on $M$
is a bilinear symmetric real-valued form $g$ 
which is $SU(2)$-invariant. Equivalently,
$g$ is quaternionic pseudo-Hermitian if
\[
g(\cdot,\cdot)= g(I\cdot, I\cdot)= g(J\cdot, J\cdot)= g(K\cdot, K\cdot).
\]
If $g$ is in addition positive definite,
$g$ is called {\bf quaternionic Hermitian}.
Notice that a quaternionic Hermitian structure exists, globally,
on any hypercomplex manifold. Indeed, one could take
any Riemannian form, and average it with $SU(2)$

As in Subsection \ref{_hc_hk_Intro_Subsection_}, 
we can associate three 2-forms 
$\omega_I$, $\omega_J$ and $\omega_K$ with $g$,
\[
\omega_I(\cdot,\cdot)= g(\cdot, I\cdot), \ \ \omega_J(\cdot,\cdot)= g(\cdot, J\cdot), \ \ \omega_K(\cdot,\cdot)= g(\cdot, K\cdot).
\]
An easy linear-algebraic calculation shows that
$\Omega_g:= \omega_J + \1 \omega_K$ has Hodge type
$(2,0)$ under $I$:
\[
  \Omega_g\in \Lambda^{2,0}_I(M).
\]
The involution $\eta \arrow J(\bar\eta)$
gives a real structure on $\Lambda^{2,0}_I(M)$.
A $(2,0)$-form $\eta$ is called {\bf real} if 
$\eta=J(\bar\eta)$.
The bundle of real $(2,0)$-forms is denoted 
$\Lambda^{2,0}_I(M,\R)$. It is easy to see that 
the form  $\Omega_g$ is real. In 
\cite{_Verbitsky_HKT-exa_}, it was shown that
the converse is also true: any real $(2,0)$-form $\eta$
is obtained from a quaternionic pseudo-Hermitian form,
which is determined uniquely from $\eta$.

\hfill

\claim\label{_pseudo-Hermi_2,0_Claim_}
Let $(M,I,J,K)$ be a hypercomplex manifold, $H$ the
bundle of quaternionic pseudo-Hermitian forms,
and $\Lambda^{2,0}_I(M,\R)$ the bundle of real
$(2,0)$-forms. Consider the map 
$H\stackrel \nu \arrow \Lambda^{2,0}_I(M,\R)$
constructed above, $\nu(g) = \Omega_g$. Then
$\nu$ is an isomorphism, and the inverse 
map is determined by $g(x, \bar y) = \Omega_g(x, J(\bar y))$,
for any $x, y \in T^{1,0}_I(M)$.

\hfill

{\bf Proof:} This is Lemma 2.10,
\cite{_Alesker_Verbitsky_HKT_}.
\endproof

\hfill

The standard diagonalization arguments, applied 
to quaternionic pseu\-do-\-Her\-mitian forms, give 
similar results about real $(2,0)$-forms
on hypercomplex manifolds.

\hfill

\proposition\label{_simu_dia_Proposition_}
Let $(M,I,J,K)$ be a hypercomplex manifold,
$\dim_\R M = 4n$, and $\eta, \eta' \in
\Lambda^{2,0}_I(M,\R)$ two real $(2,0)$-forms.
Then, locally around each point, $\eta$ and $\eta'$
can be diagonalized simultaneously: there exists
a frame 
$\xi_1, J(\bar \xi_1), \xi_2, J(\bar \xi_2), ..., \xi_n,
J(\bar \xi_n)\in \Lambda^{1,0}_I(M)$, such that
\[
\eta = \sum_i \alpha_i \xi_i \wedge J(\bar \xi_i), 
\ \ \eta' = \sum_i \beta_i \xi_i \wedge J(\bar \xi_i), 
\]
with $\alpha_i$, $\beta_i$ real-valued functions.

\hfill

{\bf Proof:} Follows from \ref{_pseudo-Hermi_2,0_Claim_}
and a standard argument which gives a simultaneous
diagonalization of two pseudo-Hermitian forms. \endproof

\hfill

In a similar spirit, 
the Gram-Schmidt orthogonalization procedure
brings the following statement.

\hfill

A real form $\eta \in \Lambda^{2,0}_I(M,\R)$
 is called {\bf strictly positive},
if it satisfies $\eta(x, J(\bar x))>0$
for any non-zero vector $x\in T^{1,0}_I(M)$.

Let $x_1, ..., x_n\in T^{1,0}_I(M)$ be a set of vector
fields. The set $\{x_i\}$ is called {\bf orthogonal with
  respect to $\eta$} if 
\[ 
\eta(x_i, x_j) = \eta(x_i, J(\bar x_j)) =0
\]
whenever $i\neq j$. 

\hfill

\proposition\label{_Gram_Schmidt_Proposition_}
(Gram-Schmidt orthogonalization procedure)
Let $\eta \in \Lambda^{2,0}_I(M,\R)$
be a real, strictly positive form on a hypercomplex
manifold, and $x_1, ..., x_n\in T^{1,0}_I(M)$
a set of vector fields, which are linearly independent everywhere.
Then there exists functions $\alpha_{i,j}$, $i>j$, such that
the vector fields 
\begin{align*}
  y_1 := & x_1,\\ 
  y_2 : = & x_2 + \alpha_{2,1} y_1, \\
  y_3 : = & x_3 + \alpha_{3,2} y_2+\alpha_{3,1} y_1, \\
  \  & ... \\
  y_k := &x_k + \sum_{i<k}\alpha_{k,i} y_i \\
  \ & ...
\end{align*}
are orthogonal.

\hfill

{\bf Proof:} Use \ref{_pseudo-Hermi_2,0_Claim_}
and apply the Gram-Schmidt orthogonalization to
the quaternionic Hermitian form associated with $\eta$.
\endproof


\section{Positive, closed $(2p,0)$-forms}


\subsection{The isomorphism 
$\Lambda^{p+q,0}_I(M) \stackrel{{\cal R}_{p,q}}\arrow\Lambda^{p,q}_{+,I}(M) $}
\label{_p+q,0_and_p,q_explicit_Subsection_}

Let $(M,I,J,K)$ be a hypercomplex manifold.
In \ref{_qD_decompo_expli_Proposition_}, an isomorphism
\[ 
  \bigoplus \Lambda^{p+q,0}_I(M)
\stackrel \Psi \arrow \bigoplus \Lambda^{p,q}_{+,I}(M)
\] 
was constructed.
As shown in \cite{_Verbitsky:HKT_}, this isomorphism is
multiplicative. It is uniquely determined by the values
it takes on $\Lambda^1(M)$: on $\Lambda^{1,0}_I(M)$,
$\Psi$ is tautological, and on $\Lambda^{0,1}_I(M)$,
we have $\Psi(x) = J(x)$. This isomorphism has an
explicit construction, which is given as follows.

\hfill

\claim\label{_R_p,q_defi_Claim_}
Let $(M,I,J,K)$ be a hypercomplex manifold, and
\[ 
  {\cal R}_{p,q}:\; \Lambda^{p+q,0}_I(M)\arrow \Lambda^{p,q}_{I}(M)
\]
map a form $\eta\in \Lambda^{p+q,0}_I(M)$
to ${\cal R}_{p,q}(\eta)$, which is defined by
\[
{\cal R}_{p,q}(\eta)(x_1, ..., x_p, \bar y_1, ..., \bar
y_q):= \eta (x_1, ..., x_p, J\bar y_1, ..., J\bar  y_q)
\]
Then ${\cal R}_{p,q}$ is multiplicative, in the following
sense:
\[
{\cal R}_{p,q}(\eta_1\wedge\eta_2) = \sum_{\begin{array}{c}{p_1+p_2=p,}\\
  {q_1+q_2=q}\end{array}} {\cal R}_{p_1,q_1}(\eta_1)\wedge {\cal R}_{p_2,q_2}(\eta_1).
\]
Moreover, ${\cal R}_{p,q}$ induces the isomorphism
\[
  \bigoplus\Lambda^{p+q,0}_I(M) 
  \stackrel{\psi}\arrow\bigoplus\Lambda^{p,q}_{+,I}(M)
\]
constructed above.

\hfill

{\bf Proof:} The multiplicativity of ${\cal R}_{p,q}$
is clear from its definition.
The isomorphism ${\cal R}$ is uniquely determined by
the values it takes on $\Lambda^1(M)$ and multiplicativity,
hence it coinsides with ${\cal R}_{p,q}$. \endproof

\hfill

This map also agrees with the differentials,
and the anticomplex involution $\eta \arrow J\bar\eta$
acting on $\Lambda^{p+q,0}_I(M)$.

\hfill

\claim \label{_R_p,q_and_6_Claim_}
Let $(M,I,J,K)$ be a hypercomplex manifold, and
\[ {\cal R}_{p,q}:\; \Lambda^{p+q,0}_I(M)\arrow \Lambda^{p,q}_{I,+}(M)\]
the map constructed in \ref{_R_p,q_defi_Claim_}.
Then
\begin{description}
\item[(i)] ${\cal R}_{p,q}(J\bar\eta)= (-1)^{pq}\overline{{\cal
    R}_{q,p}(\eta)}$.
\item[(ii)] ${\cal R}_{p,q}(\6\eta)= d^{1,0}_+{\cal  R}_{p-1,q}(\eta)$
\item[(iii)] ${\cal R}_{p,q}(\6_J\eta)= d^{0,1}_+{\cal  R}_{p,q-1}(\eta)$
\end{description}

{\bf Proof:} \ref{_R_p,q_and_6_Claim_} (i) is clear from
the definition. Using Leibniz identity, we find that
it suffices to check \ref {_R_p,q_and_6_Claim_} (ii) and (iii)
on some set of multiplicative generators of
$\bigoplus_{p,q}\Lambda^{p+q,0}_I(M)$. For functions, these identities
are clear. For $\6$-exact 1-forms,  
\ref{_R_p,q_and_6_Claim_} (ii) is clear, because
$\6^2=0$ and $(d^{1,0}_+)^2=0$, hence
\[
0 = {\cal R}_{p,q}(\6\6 f), \text{\ \ and\ \ } 
d^{1,0}_+{\cal  R}_{p-1,q}(\6f) =(d^{1,0}_+)^2f =0. 
\]
For a $\6$-exact 1-form $\eta = \6\psi$, with $\psi$
a holomorphic function, \ref{_R_p,q_and_6_Claim_} (iii)
follows from
\[
{\cal R}_{p,q}(\6_J\6\psi) = - {\cal R}_{p,q}(\6\6_J\psi)
= - {\cal R}_{p,q}(\6J\bar\6\psi)=0.
\]
The functions, together with 1-forms $\eta = \6\psi$, with $\psi$
a holomorphic function, generate the algebra $\Lambda^{*,0}_I(M)$ 
multiplicatively. Now, the Leibniz identity can be used to prove that
\ref {_R_p,q_and_6_Claim_} (ii) and (iii) is true on the whole
$\Lambda^{*,0}_I(M)$.

Please notice that we just gave a proof of \ref{_bico_ide_Theorem_}.
\endproof

\subsection{Strongly positive, weakly positive and real $(2p,0)$-forms}
\label{_positive_2p,0_Subsection_}

The notion of positive $(2p,0)$-forms on hypercomplex
manifolds was developed
in \cite{_Alesker_Verbitsky_HKT_} and in ongoing
collaboration with S. Alesker. 

\hfill

Let $\eta\in \Lambda^{p,q}_I(M)$ be a 
differential form. Since $I$ and $J$ anticommute,
$J(\eta)$ lies in $\Lambda^{q,p}_I(M)$.
Clearly, $J^2\restrict {\Lambda^{p,q}_I(M)}=(-1)^{p+q}$.
For $p+q$ even, $J\restrict {\Lambda^{p,q}_I(M)}$
is an anticomplex involution, that is, a real structure
on $\Lambda^{p,q}_I(M)$. 
A form 
$\eta \in \Lambda^{2p,0}_I(M)$ is called {\bf real} if
$J(\bar\eta)=\eta$. We denote the bundle of
real $(2p,0)$-forms by $\Lambda^{2p,0}_I(M, \R)$.

For a real $(2p,0)$-form, 
\begin{multline}\label{_real_forms_real_Equation_}
   \eta\left(x_1, J(\bar x_1), x_2, J(\bar x_2), ... x_p,
   J(\bar x_p)\right)=\\ = 
   \bar \eta\left(J(x_1), J^2 (\bar x_1), J(x_2), J^2(\bar x_2), ... J(x_p),
   J^2(\bar x_p)\right)=\\ = 
 \bar \eta\left(\bar x_1, J(x_1), \bar x_2, J(x_2), ...\bar  x_p,
   J(x_p)\right),
\end{multline}
for any $x_1, ..., x_p \in T^{1,0}_I(M)$.
From \eqref{_real_forms_real_Equation_},
we obtain that the number 
\[ \eta\left(x_1, J(\bar x_1), x_2, J(\bar x_2), ... x_p,
   J(\bar x_p)\right)
\] 
is always real.

\hfill

\definition
Let $(M,I,J,K)$ be a hypercomplex manifold, and
$\eta \in \Lambda^{2p,0}_I(M)$ a real $(2p,0)$-form.
It is called {\bf weakly positive}, if 
\[ \eta(x_1, J(\bar x_1), x_2, J(\bar x_2), ..., x_p,
   J(\bar x_p)) \geq 0,
\]
for any $x_1, ..., x_p \in T^{1,0}_I(M)$.

\hfill

Let $\dim_\R M=4n$. The complex line bundle
$\Lambda^{2n,0}(M)$ is equipped with a real structure,
hence it is a complexification of a real line bundle
$\Lambda^{2n,0}_I(M, \R)$. This real line bundle
is trivial topologically. To see this, take
a quaternionic Hermitian form $q$ on $M$
(such a form always exists: see Section
\ref{_quat_pseudo-H_Section_}).
Let $\Omega:= \omega_J + \1\omega_K$ be
the corresponding $(2,0)$-form. Since
$J\omega_J = \omega_J$, $J(\omega_K)=-\omega_K$,
the form $\Omega$ is real. Then, $\Omega^n$
is a nowhere degenerate, real section
which trivializes $\Lambda^{2n,0}_I(M, \R)$.

The pairing 
\[
\Lambda^{2p,0}_I(M,\R)\times \Lambda^{2n-2p,0}_I(M,\R) \arrow
\Lambda^{2n,0}_\R(M,\R)
\]
is nowhere degenerate. Denote by ${\cal
  C}_w\subset\Lambda^{2*,0}_I(M,\R)$  the cone of
weakly positive forms, and ${\cal
  C}_s\subset\Lambda^{2*,0}_I(M,\R)$ 
the dual cone. This cone is called
{\bf the cone of strongly positive forms}.

This notion is well known in complex geometry;
a complex analogue of the following claim
is often used as a definition of strongly
positive cone, and then the above definition
becomes a (trivial) theorem.

\hfill

\claim
Let $M$ be a hypercomplex manifold. The 
cone ${\cal C}_s\subset\Lambda^{2*,0}_I(M,\R)$ 
of strongly positive real $(2p,0)$-forms
is multiplicatively generated by products of forms $\xi\wedge J(\bar\xi)$, for
$\xi\in \Lambda^{1,0}_I(M)$. 

\hfill

{\bf Proof:} A form $\eta$ is weakly positive if
\[ 
 \langle\eta, \xi_1\wedge J(\bar \xi_1)\wedge \xi_2\wedge J(\bar
\xi_2)\wedge  ... \wedge J(\bar \xi_p)\rangle \geq 0
\]
for any $\xi_1, ..., \xi_p\in \Lambda^{1,0}_I(M)$.
Therefore, weakly positive cone is dual to the cone
generated by such products.
\endproof

\hfill

The strong positivity of a form implies its weak positivity.
Unlike the complex case, in the quaternionic case this is
not immediate from its definition.

For $p=n$, this implication can be seen as follows.
For any $\xi_1, ..., \xi_p\in \Lambda^{1,0}_I(M)$, we have
\[ \xi_1\wedge J(\bar \xi_1)\wedge \xi_2\wedge J(\bar
\xi_2)\wedge  ... \wedge J(\bar \xi_n)= \frac 1 {n!}\Omega^n,
\]
where $\Omega = \sum \xi_i\wedge J(\bar \xi_1)$
is a $(2,0)$-form, which is obtained from 
a quaternionic Hermitian form $q$ as in \ref{_pseudo-Hermi_2,0_Claim_}.
The form $\Omega^n$ is positive, because 
for $\{\langle x_i, J(\bar x_i)\}$ 
pairwise orthogonal with respect to $q$, we have  
\[
\Omega^n(x_1, J(\bar x_1), ..., x_n, J(\bar x_n))=
\prod_i q(x_i, \bar x_i),
\]
and for $\{x_i\}$ non-orthogonal, this set can be orthogonalized,
without changing $\eta(x_1, J(\bar x_1), ..., x_n, J(\bar x_n))$,
as shown in \ref{_Gram_Schmidt_Proposition_}. 

This gives
\begin{equation}\label{_Omega^n_positive_Equation_}
\frac 1 {n!}\Omega^n(x_1, J(\bar x_1), ..., x_n, J(\bar x_n)\geq 0
\end{equation}

For $p<n$, we restrict $\eta$ to a quaternionic
subspace generated by $x_1, ... x_p$, and
find that the positivity of 
\[
\xi_1\wedge J(\bar \xi_1)\wedge \xi_2\wedge J(\bar
\xi_2)\wedge  ... \wedge J(\bar \xi_p)
\left(x_1, J(\bar x_1), x_2, J(\bar x_2), ... x_p,
   J(\bar x_p)\right)
\]
follows from \eqref{_Omega^n_positive_Equation_}.

\hfill

Recall that a real $(p,p)$-form $\rho$
on a complex manifold $X$ is called {\bf weakly positive}
if 
\[ (-\1)^p\rho(x_1, \bar x_1, ... x_p, \bar x_p)\geq 0,
\]
for any $x_1, ... x_p \in T^{1,0}(X)$.

\hfill

\claim\label{_posi_real_R_p,p_Claim_}
Let $(M,I,J,K)$ be a hypercomplex manifold,
and \[ {\cal R}_{p,p}:\; \Lambda^{2p,0}_I(M) \arrow
\Lambda^{p,p}_I(M)\] the map constructed in 
Subsection \ref{_p+q,0_and_p,q_explicit_Subsection_}.
Consider a $(2p,0)$-form \\ $\eta \in \Lambda^{2p,0}_I(M)$.
Then
\begin{description}
\item[(i)] $\eta$ is real if and only if $(\1)^p{\cal R}_{p,p}(\eta)$
is real (in the usual sense).
\item[(ii)] $\eta$ is weakly positive if and only if
$(\1)^p{\cal R}_{p,p}(\eta)$ is a weakly positive $(p,p)$-form. 
\end{description}

{\bf Proof:} \ref{_posi_real_R_p,p_Claim_} (i) is clear from 
the definition.  Indeed,
\[ 
  {\cal R}_{p,p}(\eta)(x_1,  \bar x_1, ..., x_p,  \bar x_p)
= \eta(x_1,  J(\bar x_1), ..., x_p,  J(\bar x_p)).
\]
It is easy to see that a $(p,p)$-form $\rho$
is real if and only if $(\1)^p\rho$
satisfies $\rho(x_1,  \bar x_1, ..., x_p,  \bar x_p)\in \R$.

\ref{_posi_real_R_p,p_Claim_} (ii) is also clear.
Indeed,
\begin{multline*}
\eta\left(x_1, J(\bar x_1), x_2, J(\bar x_2), ... x_p,
   J(\bar x_p)\right) = \\ = (-1)^{p(p-1)}\eta\left(x_1, x_2, ..., x_p, J(\bar x_1), 
  J(\bar x_2), ..., J(\bar x_p)\right). 
\end{multline*}
Therefore,
\begin{multline}\label{_R_p,p_positive_Equation_}
  {\cal R}_{p,p}(\eta)(x_1,  \bar x_1, ..., x_p,  \bar x_p)
{\cal R}_{p,p}(\eta)(x_1, ..., x_p, \bar x_1, ..., \bar x_p)= \\ =
  \eta\left(x_1, ..., x_p, J(\bar x_1), ..., J(\bar x_p)\right)=
  \eta\left(x_1, J(\bar x_1), x_2, J(\bar x_2), ... x_p,
   J(\bar x_p)\right)
\end{multline}
Then, \eqref{_R_p,p_positive_Equation_} is non-negative if and
only if $\eta$ is weakly positive, and this is equivalent to 
$(\1)^p{\cal R}_{p,p}(\eta)$ being weakly positive, by definition
of positive $(p,p)$-forms.
\endproof

\subsection{The map ${\cal V}_{p,q}:\;
  \Lambda^{p+q,0}_I(M)\arrow\Lambda^{n+p, n+q}_I(M)$\\
on $SL(n, {\Bbb H})$-manifolds}
\label{_V_p,q_Subsection_}

Let $(M,I,J,K)$ be a hypercomplex manifold, $\dim_\R M =4n$,
and 
\[ 
  {\cal R}_{p,q}:\; \Lambda^{p+q,0}_I(M)\arrow \Lambda^{p,q}_{I,+}(M)
\]
the isomorphism defined in Subsection
\ref{_p+q,0_and_p,q_explicit_Subsection_}.
Consider the projection 
\begin{equation}\label{_proj_to_+_Equation_} 
\Lambda^{p,q}_{I}(M)\arrow
\Lambda^{p,q}_{I,+}(M),
\end{equation}
and let
\[ 
  R:\; \Lambda^{p,q}_{I}(M)\arrow\Lambda^{p+q,0}_I(M)
\]
denote the composition of \eqref{_proj_to_+_Equation_} 
and ${\cal R}_{p,q}^{-1}$. 

\hfill

\lemma\label{_R_expli_and_proj_Lemma_}
In these assumptions,
\begin{equation}\label{_R_expli_Equation_}
   R(\xi_1\wedge ... \wedge \xi_p \wedge \bar\xi_{p+1}
   \wedge ... \wedge \bar\xi_{p+q}) = 
   \xi_1  ... \wedge \xi_p \wedge J(\bar\xi_{p+1})
   \wedge ... \wedge J(\bar\xi_{p+q}),
\end{equation}
for any $\xi_1, ..., \xi_{p+q} \in \Lambda^{1,0}_I(M)$.

\hfill

{\bf Proof:}
Denote by $R'$ the map defined by the formula
\eqref{_R_expli_Equation_}. From the definition of the
$SU(2)$-action on $\Lambda^*(M)$ it is apparent that
$R'(\eta)$ belongs to the same $SU(2)$-representation
as $\eta$. Since $R'(\eta)$ lies in
$\Lambda^{p+q,0}_I(M)$, it belongs to
$\Lambda^*_+(M)$. Therefore, $R'$ vanishes on the
kernel of \eqref{_proj_to_+_Equation_}. By definition,
$R$ is the unique map $\Lambda^{p,q}_{I}(M)\arrow\Lambda^{p+q,0}_I(M)$
vanishing on the kernel of \eqref{_proj_to_+_Equation_} 
and satisfying 
\[  R\circ {\cal R}_{p,q}= \Id_{\Lambda^{p+q,0}_I(M)}.\]
To prove that $R'=R$ it suffices now to check that
$R({\cal R}_{p,q}(\eta))=\eta$, but this
is obvious from the definition.
\endproof

\hfill

\remark
The formula \eqref{_R_expli_Equation_}
could be used as a definition of $R$.

\hfill

The map $R$ is compatible with Dolbeault differentials,
in the following sense.

\hfill

\lemma\label{_R_differe_Lemma_}
Let $(M,I,J,K)$ be a hypercomplex manifold,
and 
\[ 
  R:\; \Lambda^{p,q}_{I}(M)\arrow\Lambda^{p+q,0}_I(M)
\]
the map defined above. Then
\begin{equation}\label{_R_differe_Equation_}
R(\6\eta) = \6R(\eta), \text{\ \ and \ \ }
R(\bar\6\eta) = \6_J R(\eta).
\end{equation}
{\bf Proof:}
\ref {_R_differe_Lemma_}
follows immediately from \ref{_R_p,q_and_6_Claim_}
and $R\circ {\cal R}_{p,q}= \Id_{\Lambda^{p+q,0}_I(M)}$,
which is a part of the definition of $R$.
\endproof

\hfill

Let $\Phi_I$ be a nowhere degenerate 
holomorphic section of $\Lambda^{2n,0}_I(M)$. Assume that $\Phi_I$ is
real, that is, $J(\Phi_I)=\bar\Phi_I$, and positive.

Existence of such a section is highly non-trivial.
When $M$ is hyperk\"ahler, we could take the
top power of the holomorphic symplectic form
$\Omega=\omega_J+\1\omega_K$. For a general
hypercomplex $M$, such a form $\Phi_I$ is
preserved by the Obata connection, and reduces
the holonomy of Obata connection to a subgroup
of $SL(n,{\Bbb H})$. Such manifolds were studied
in \cite{_Verbitsky:canoni_} and \cite{_BDV:nilmanifolds_}.

A manifold with a nowhere degenerate, real, positive
form $\Phi_I\in \Lambda^{2n,0}_I(M)$ is called
{\bf an $SL(n, {\Bbb H})$-manifold}.

\hfill

\remark
Let $(M,I,J,K, \Phi_I)$ be an $SL(n, {\Bbb H})$-manifold.
For any section $\eta \in\Lambda^{2n,0}_I(M)$,
positivity of $\eta$ in the quaternionic sense is equivalent
to positivity of $\eta \wedge \Phi_I\in \Lambda^{2n,2n}_I(M)$, 
in the usual sense.

\hfill

Define the map
\[ {\cal V}_{p,q}:\;
  \Lambda^{p+q,0}_I(M)\arrow\Lambda^{n+p, n+q}_I(M)
\]
by the relation
\begin{equation}\label{_V_p,q_via_test_form_Equation_}
{\cal V}_{p,q}(\eta) \wedge \alpha = \eta \wedge R(\alpha)\wedge \bar\Phi_I,
\end{equation}
for any test form $\alpha \in \Lambda^{n-p, n-q}_I(M)$.

\hfill

The map ${\cal V}_{p,p}$ is especially remarkable,
because it maps closed, positive
$(2p,0)$-forms to closed, positive $(n+p, n+p)$-forms,
as the following proposition implies.

\hfill

\proposition\label{_V_main_Proposition_}
Let $(M,I,J,K, \Phi_I)$ be an $SL(n, {\Bbb H})$-manifold, and
\[ {\cal V}_{p,q}:\;
  \Lambda^{p+q,0}_I(M)\arrow\Lambda^{4n-p, 4n-q}_I(M)
\]
be the map defined above.
Then
\begin{description}
\item[(i)] ${\cal V}_{p,q}(\eta)= {\cal R}_{p,q}(\eta) \wedge {\cal V}_{0,0}(1)$.
\item[(ii)]  The map ${\cal V}_{p,q}$ is injective, for
  all $p$, $q$.
\item[(iii)] $(\1)^{(n-p)^2}{\cal V}_{p,p}(\eta)$ is real if and
  only $\eta\in\Lambda^{2p,0}_I(M)$ is real, 
and weakly positive if and only if $\eta$ is weakly positive.
\item[(iv)] ${\cal V}_{p,q}(\6\eta)= \6{\cal V}_{p-1,q}(\eta)$,
and ${\cal V}_{p,q}(\6_J\eta)= \bar\6{\cal  V}_{p,q-1}(\eta)$.
\item[(v)] ${\cal V}_{0,0}(1) = \lambda {\cal
  R}_{n,n}(\Phi_I)$, where $\lambda$ is a positive rational number,
depending only on the dimension $n$.
\end{description}

{\bf Proof:} The map 
$R:\; \Lambda^{p,q}_{I}(M)\arrow\Lambda^{p+q,0}_I(M)$
is by construction multiplicative, and satisfies
\begin{equation}\label{_R_R_p,q_inverse_Equation_}
R({\cal R}_{p,q}(\eta))=\eta,
\end{equation}
for all $\eta\in\Lambda^{p+q,0}_I(M)$.
This gives
\begin{equation}\label{_V_multi_calc_Equation_}
{\cal V}_{p,q}(\eta) \wedge\alpha = \eta \wedge R(\alpha) \wedge
\Phi_I= R({\cal R}_{p,q}(\eta)\wedge\alpha) \wedge
\Phi_I = {\cal V}_{0,0}(1)\wedge{\cal R}_{p,q}(\eta)\wedge\alpha
\end{equation}
(to obtain the last equation, we take the test-form
$\alpha':={\cal R}_{p,q}(\eta)\wedge\alpha$ and apply
\eqref{_V_p,q_via_test_form_Equation_}).
Since $\alpha$ is arbitrary, \eqref{_V_multi_calc_Equation_}
gives \[ {\cal V}_{p,q}(\eta)={\cal V}_{0,0}(1)\wedge{\cal R}_{p,q}(\eta).\]
This proves \ref{_V_main_Proposition_} (i).

Injectivity of ${\cal V}_{p,q}$ is clear, because for any
$\eta\in\Lambda^{p+q,0}_I(M)$ there exists $\chi$ such that 
$\eta \wedge \chi \wedge \Phi_I\neq 0$.
Using \eqref{_R_R_p,q_inverse_Equation_}, we find that
\[ 
  {\cal V}_{p,q}(\eta) \wedge {\cal R}_{n-p,n-q}(\chi)=
  \eta \wedge R({\cal R}_{n-p,n-q}(\chi)) \wedge \Phi_I=
  \eta \wedge \chi \wedge \Phi_I\neq 0.
\]
We proved \ref{_V_main_Proposition_} (ii).

From \ref{_R_p,q_and_6_Claim_} (i),
we obtain that $R(\bar\alpha)= (-1)^{pq} R(\alpha)$,
for any $\alpha\in \Lambda^{p,q}_I(M)$. Then
\[ 
  {\cal V}_{p,q}(J\bar\eta)=
  (-1)^{(n-p)(n-q)}\overline{{\cal V}_{q,p}(\eta)}
\]
as follows from \eqref{_V_p,q_via_test_form_Equation_}.
Then, $(\1)^{p}{\cal V}_{p,p}(\eta)$ is real if
 $J\bar\eta=\eta$. The ``only if'' part follows from
injectivity of ${\cal V}_{p,p}$.

To check the weak positivity of $(\1)^{p}{\cal V}_{p,p}$,
take  $\alpha = \xi_1 \wedge \bar \xi_1 \wedge ... \wedge
\xi_{n-p}\wedge \xi_{n-p}$, with $\xi_1, ..., \xi_{n-p} \in
\Lambda^{1,0}_I(M)$. Then $(-\1)^{n-p}\alpha$ is positive.
Such forms generate the strongly positive cone. Then
$R(\alpha)= \xi_1 \wedge J(\bar \xi_1) \wedge ... \wedge
\xi_{n-p}\wedge J(\bar \xi_{n-p})$ is strongly positive
by definition, and, moreover, $R(\alpha)$, for all such $\alpha$,
 generate the strongly positive cone.

The weak positivity of 
$(-\1)^{n-p}{\cal V}_{p,q}(\eta)$ is equivalent to
\[
(-\1)^{n-p}{\cal V}_{p,q}(\eta) \wedge \alpha \geq 0,
\]
and the weak positivity of $\eta$ 
is equivalent to 
\[ 
  \eta \wedge R(\alpha)\wedge \bar\Phi_I \geq 0.
\]
These two inequalities are equivalent by the formula
\eqref{_V_p,q_via_test_form_Equation_} which is a definition of
${\cal V}_{p,q}(\eta)$. We proved \ref{_V_main_Proposition_} (iii).

\ref{_V_main_Proposition_} (iv) follows from the 
Stokes' formula
\[ 
  \int_M \6\alpha\wedge \beta = (-1)^{\deg \alpha}\int_M \alpha\wedge \6\beta,
\]
where $\alpha$ or $\beta$ have compact support.

Take an $(n-q, n-p)$-form 
$\alpha$ with compact support. By \ref{_R_differe_Lemma_},
\begin{align*}
\int_M
 {\cal V}_{p,q}(\6\eta) \wedge \alpha =& 
 \int_M\6\eta \wedge R(\alpha)\wedge \bar\Phi_I =
 (-1)^{p+q-1} \int_M\eta \wedge \6R(\alpha)\wedge \bar\Phi_I =\\ = &
 (-1)^{p+q-1} \int_M\eta \wedge R(\6\alpha)\wedge \bar\Phi_I=\\ = &
  (-1)^{p+q-1} \int_M{\cal V}_{p-1,q}(\eta) \wedge \6\alpha =\\ = &
  \int_M
 \6{\cal V}_{p-1,q}(\eta) \wedge \alpha.
\end{align*}
Applying complex conjugation to both sides of
${\cal V}_{p,q}(\6\eta)= \6{\cal V}_{p-1,q}(\eta)$
and using 
\[ 
  {\cal V}_{p,q}(J\bar\eta)=
  (-1)^{(n-p)(n-q)}\overline{{\cal V}_{q,p}(\eta)}
\]
and $J\bar\6\eta = \6_JJ(\bar\eta)$,
we obtain the second equation of \ref{_V_main_Proposition_} (iv).

\ref{_V_main_Proposition_} (v) follows from a direct (but tedious)
linear-algebraic
calculation. The bundle $\Lambda^{n,n}_{I,+}(M)$ is 1-dimensional,
by \ref{_qD_decompo_expli_Proposition_}. The form
${\cal V}_{0,0}(1)$ lies in $\Lambda^{n,n}_{I,+}(M)$. Indeed,
\[
{\cal V}_{0,0}(1)\wedge \alpha = R(\alpha) \wedge \bar\Phi_I,
\]
and therefore $\alpha\arrow {\cal V}_{0,0}(1)\wedge \alpha$
vanishes on all forms of weight less than $2n$. Therefore,
${\cal V}_{0,0}(1)$ has weight $2n$, hence belongs to 
$\Lambda^{n,n}_{I,+}(M)$. The form ${\cal
  R}_{n,n}(\Phi_I)$ is a nowhere degenerate section
of $\Lambda^{n,n}_{I,+}(M)$, by construction; therefore,
${\cal V}_{0,0}(1)$ is proportional to ${\cal
  R}_{n,n}(\Phi_I)$:
\[
{\cal V}_{0,0}(1)= \lambda {\cal
  R}_{n,n}(\Phi_I),
\]
where $\lambda$ is a smooth function on $M$. 
To prove \ref{_V_main_Proposition_} (v),
we need to show that $\lambda$ is a positive rational
number depending only from $n$. 
Since $(\1)^n{\cal R}_{n,n}(\Phi_I)$ 
and $(\1)^n{\cal V}_{0,0}(1)$ are both real and positive,
by \ref{_V_main_Proposition_} (iii) and 
\ref{_posi_real_R_p,p_Claim_}, $\lambda$ 
is real and positive. Taking $\alpha=\Phi_I$ and aplying 
\eqref{_V_p,q_via_test_form_Equation_}, we obtain
\[
1\wedge \Phi_I \wedge \bar \Phi_I = R({\cal
  R}_{n,n}(\Phi_I)) \wedge \bar \Phi_I = {\cal V}_{0,0}(1)
\wedge {\cal R}_{n,n}(\Phi_I) = \lambda 
{\cal R}_{n,n}(\Phi_I))\wedge {\cal R}_{n,n}(\Phi_I)
\]
This gives an expression for $\lambda$:
\[
\lambda = \frac{\Phi_I \wedge \bar \Phi_I}
{{\cal R}_{n,n}(\Phi_I)\wedge {\cal R}_{n,n}(\Phi_I)}.
\]
From this formula, it is clear that $\lambda$ is
independent from the choice of $\Phi_I$.
Therefore, we may assume that $\Phi_I$ is associated
with a quaternionic Hermitian form $q$ as above:
$\Phi_I= \Omega^n$, where $\Omega= \omega_J+ \1\omega_K$,
and $\omega_J, \omega_K$ are the Hermitian skew-linear 
forms of $(M,J)$ and $(M,K)$. From
the definition of ${\cal R}_{p,q}$, it is clear that
${\cal R}_{1,1}(\Omega)=\omega_I$. Using multiplicativity of
${\cal R}_{p,p}$, we obtain
\[
{\cal R}_{n,n}(\Omega^n) = \Pi_+({\cal R}_{1,1}(\Omega)^n) =
\Pi_+(\omega_I^n),
\]
where $\Pi_+$ is the $SU(2)$-invariant 
projection to the $\Lambda^*_+(M)$-part. Since
the metric on $\Lambda^*(M)$ is $SU(2)$-invariant,
the weight decomposition of $\Lambda^*(M)$ is orthogonal;
therefore, $\Pi_+$ is an orthogonal projection to
$\Lambda^*_+(M)$.

Consider the algebra $A^*= \oplus A^{2i}$
generated by $\omega_I$, $\omega_J$, and $\omega_K$.
In \cite{_Verbitsky:Symplectic_II_}, this algebra 
was computed explicitly. It was shown, that, up to the middle
degree, $A^*$ is a symmetric algebra with generators
$\omega_I$, $\omega_J$, $\omega_K$. The algebra $A^*$ has Hodge
bigrading $A^k = \bigoplus\limits_{p+q=k}A^{p,q}$, and its
$A^{p,p}$-part is generated by the forms
\[ \omega_I^i\wedge (\Omega\wedge\bar\Omega)^j,
\] 
$i, j = 0, 1, 2, ...$
From the Clebsch-Gordan formula, we obtain that
$A^{2i}_+:= \Lambda^{2i}_+(M)\cap A^{2i}$, for $i\leq n$,
is an orthogonal complement to $Q(A^{2i-4})$,
where $Q(\eta) = \eta \wedge (\omega_I^2 + \omega_J^2+\omega_K^2)$.
The space $A^{n,n}_+= \ker Q^* \restrict {A^{n,n}}$ is 1-dimensional,
as we have shown above, and generated by 
${\cal R}_{n,n}(\Omega^n)$. Clearly,
\[
Q^*\left(\omega_I^i\wedge (\Omega\wedge\bar\Omega)^j\right) =
\omega_I^{i-2}\wedge (\Omega\wedge\bar\Omega)^j +
\omega_I^i\wedge (\Omega\wedge\bar\Omega)^{j-2}.
\]
Therefore, $\ker Q^* \restrict A^{n,n}$ is generated by
\begin{equation}\label{_Xi_expli_Equation_}
\Xi:= \omega_I^n - \omega_I^{n-2}\wedge (\Omega\wedge\bar\Omega) +
\omega_I^{n-4}\wedge (\Omega\wedge\bar\Omega)^2-
\omega_I^{n-6}\wedge (\Omega\wedge\bar\Omega)^3 + ... 
\end{equation}
Since ${\cal R}_{n,n}(\Omega^n)$ is equal to the
projection of
$\omega_I^n$ to $\ker Q^*$, this gives
\[
{\cal R}_{n,n}(\Omega^n)=\Xi\cdot \frac{(\omega_I^n,
  \Xi)}{(\Xi,\Xi)}=\gamma\Xi,
\]
where $\gamma$ is a rational coefficient which can
be expressed through binomial
coefficients using \eqref{_Xi_expli_Equation_}. 
A similar calculation can be used
to express 
\[ 
\lambda = \frac{\Phi_I \wedge \bar \Phi_I}
{{\cal R}_{n,n}(\Phi_I))\wedge {\cal R}_{n,n}(\Phi_I)} =
\frac{\Omega^n \wedge \bar\Omega^n}{\gamma^2 \Xi\wedge\Xi}
\]
through a combinatorial expression which would take half a page.
\endproof


\section{Sibony's Lemma for positive $(2p,0)$-forms}
\label{_Sibony_Section_}


\subsection{$\omega^q$-positive (1,1)-forms}

Recall that a real $(p,p)$-form $\eta$
on a complex manifold is called {\bf weakly positive}
if for any complex subspace $V\subset T_c M$, 
$\dim_\C V=p$, the restriction $\rho\restrict V$
is a non-negative volume form. Equivalently,
this means that 
\[ 
  (\1)^p\rho(x_1, \bar x_1, x_2, \bar x_2, ... x_p, \bar
  x_p)\geq 0,
\]
for any vectors $x_1, ... x_p\in T_x^{1,0}M$.
A form is called {\bf strongly positive} if it can 
be expressed as a sum
\[
\eta = (\1)^p\sum_{i_1, ... i_p} 
\alpha_{i_1, ... i_p} \xi_{i_1} \wedge \bar\xi_{i_1}\wedge ... 
\wedge \xi_{i_p} \wedge \bar\xi_{i_p}, \ \  
\]
running over some set of $p$-tuples 
$\xi_{i_1}, \xi_{i_2}, ..., \xi_{i_p}\in \Lambda^{1,0}(M)$,
with $\alpha_{i_1, ... i_p}$ real and non-negative functions on $M$.

The strongly positive and the weakly positive forms
form closed, convex cones in the space 
$\Lambda^{p,p}(M,\R)$ of real $(p,p)$-forms.
These two cones are dual with respect to the Poincare pairing
\[
\Lambda^{p,p}(M,\R) \times \Lambda^{n-p,n-p}(M,\R)\arrow \Lambda^{n,n}(M,\R)
\]
where $n=\dim_\C M$.
For (1,1)-forms and $(n-1,n-1)$-forms,
the strong positivity is equivalent
to weak positivity.

\hfill

\definition
Let $(M, \omega)$ be a K\"ahler manifold.
A real (1,1)-form $\eta\in \Lambda^{1,1}(M, \R)$
is called $\omega^q$-positive if $\omega^{q-1}\wedge \eta$
is a weakly positive form. 

\hfill

This notion was studied in \cite{_Verbitsky:omega-psh_},
in connection with plurisubharmonic functions on
calibrated manifolds (\cite{_Harvey_Lawson:Psh_},
\cite{_Harvey_Lawson:Dua_}). 
In \cite{_Verbitsky:omega-psh_}, a characterization of
$\omega^q$-positivity in terms of the eigenvalues was
obtained.  At each point $x\in M$, we can
find an orthonormal basis $ \xi_1, ... \xi_n \in
\Lambda^{1,0}_x(M)$,  such that
\[
\eta = -\1\sum_i \alpha_i \xi_i\wedge\bar\xi_i.
\]
The numbers $\alpha_i$ are called {\bf the eigenvalues}
of $\eta$ at $x$.

\hfill

The following theorem was proven in \cite{_Verbitsky:omega-psh_}.

\hfill

\theorem\label{_omega^q_psh_via_eigenva_Theorem_}
Let $(M, \omega)$ be a K\"ahler manifold, and
$\eta\in \Lambda^{1,1}(M, \R)$ a real (1,1)-form.
 Let $\alpha_1(x), \alpha_2(x), ..., \alpha_n(x)$ denote the
eigenvalues of $\eta$ at $x\in M$.
Then the following conditions are equivalent.
\begin{description}
\item[(i)] $\eta$ is $\omega^q$-positive
\item[(ii)] $\eta\wedge \omega^{q-1}$ is weakly
positive
\item[(ii)] $\eta\wedge \omega^{q-1}$ is strongly
positive
\item[(iv)] The sum of any $q$ eigenvalues of $\eta$ is
  positive, for any $x\in M$:
\begin{equation}\label{_sum_smallest_Equation_}
\sum_{k=1}^q \alpha_{i_k}(x) \geq 0,
\end{equation}
for any  
$q$-tuple $\{i_1, ... i_q\} \subset \{ 1, 2, ... , n\}$.
\end{description}

{\bf Proof:} This is \cite{_Verbitsky:omega-psh_}, Theorem 2.4. In \cite{_Verbitsky:omega-psh_}, 
this statement was stated for forms $\eta = dd^c\phi$, 
but the proof is purely linear-algebraic, and can be 
extended to arbitrary (1,1)-forms. \endproof

\hfill

\definition
A form $\eta$ is called {\bf strictly  $\omega^q$-positive},
if $\eta - h\omega$ is $\omega^q$-positive, for some
continuous, nowhere vanishing, positive function $h$ on $M$.

\subsection{Positive $(2p,0)$-forms on hypercomplex
  manifolds}

Let $(M, I, J, K)$ be a hypercomplex manifold.
In Subsection
\ref{_positive_2p,0_Subsection_}, a notion of positivity for $(2p, 0)$-forms
on $M$ was defined. 
We say that
a real $(2,0)$-form $\eta$ is $\Omega^q$-positive if
$\eta \wedge \Omega^{q-1}$ is positive, and {\bf strictly 
positive} if $\eta \wedge \Omega^{q-1}- h \Omega^q$
is positive, for some continuous, nowhere vanishing, 
positive function $h$ on $M$.

As shown in \ref{_pseudo-Hermi_2,0_Claim_}, quaternionic
pseudo-Hermitian forms are in (1,1)-correspondence
with real $(2,0)$-forms. This allows one
to diagonalize a given $(2,0)$-form $\eta$
locally in an orthonormal frame
(\ref{_simu_dia_Proposition_}).

Given a real $(2,0)$-form $\eta$ on a hyperk\"ahler
manifold, at any point $x\in M$
there exists an orthonormal frame
$\xi_1, J\bar\xi_1, ..., \xi_n, J\bar\xi_n\in \Lambda^{1,0}_I(M)$,
such that $\eta\restrict x$ is written as
\[
\eta\restrict x = \sum_i\alpha_i \xi_1\wedge J\bar\xi_1,
\]
with $\alpha_i$ being real-valued functions.
The condition of $\Omega^q$-positivity is equivalent
to the inequality 
\begin{equation}\label{_sum_smallest_Omega_Equation_}
\sum_{k=1}^q \alpha_{i_k}(x) \geq 0,
\end{equation}
just like in \ref{_omega^q_psh_via_eigenva_Theorem_}.

Given a (1,1)-form  $\eta\in\Lambda^{1,1}_I(M)$,
consider a $(2,0)$-form $R(\eta)\in
\Lambda^{2,0}_I(M)$, 
\[ R(\eta)(x, y) := \eta(x, J(y)).\]
Clearly, $R(\eta)$ is real and positive if
$\eta$ is real and positive.
It is easy to see that $R$ vanishes on
$SU(2)$-invariant forms, and induces
an isomorphism $\Lambda^{1,1}_{+,I}(M)\arrow \Lambda^{2,0}_I(M)$
described in \ref{_R_p,q_defi_Claim_}
(see \ref{_R_expli_and_proj_Lemma_} for a detailed argument).

\hfill

\lemma\label{_R_eta_Omega_posi_Lemma_}
Let $M$ be a hyperk\"ahler manifold, $\dim_\R M =4n$, and
$\eta \in \Lambda^{1,1}(M,\R)$ a real $(1,1)$-form,
which is $\omega^{2n-2p}$-positive. Then
$R(\eta)$ is $\Omega^{n-p}$-positive.

\hfill

{\bf Proof:} Denote by $\eta'$ the $(1,1)$-form
$\eta - \eta_{inv}$, where $\eta_{inv}= \frac 1 2 (\eta + J(\eta))$
denotes the $SU(2)$-invariant part of $\eta$. Clearly,
\[
\eta' = \frac 1 2 (\eta - J(\eta)).
\]
Since $-J(\eta)$ has the same eigenvalues as $\eta$,
by \ref{_omega^q_psh_via_eigenva_Theorem_} (iv)
it is also $\omega^{2n-2p}$-positive. Then
$\eta'$ is $\omega^{2n-2p}$-positive, too.

Using the orthonormal frame as
in the proof of \eqref{_sum_smallest_Omega_Equation_},
we find that $\eta'$ can be written
as 
\[
\eta' =-\1\sum_i \alpha_i \xi_i\wedge\bar\xi_i,
\]
with $\xi_i$ an orthonormal basis in $\Lambda^{1,0}_I(M)$
satisfying 
\[ J(\xi_{2i-1})= \bar\xi_{2i}, \ \ J(\xi_{2i})=
   -\bar\xi_{2i-1}
\]
(see \ref{_simu_dia_Proposition_}). 
Since $J(\eta')=-\eta'$, the eigenvalues of $\eta'$
occur in pairs: 
\begin{equation}\label{_in_pairs_Equation_}
\alpha_{2i-1}=\alpha_{2i}. 
\end{equation}
Renumbering
the basis, we may assume that 
$\alpha_1 \leq \alpha_2\leq ... \leq \alpha_{2n}$.
Now, $\omega^{2n-2p}$-positivity of $\eta'$ is equivalent
to 
\begin{equation}\label{_first_2n-2k_eigen_Equation_}
\alpha_1 + \alpha_2 + ... + \alpha_{2n-2p}\geq 0.
\end{equation}
By definition,
\[
R(\eta')= 2 \sum_i \alpha_{2i} \xi_{2i-1}\wedge \xi_{2i},
\]
hence \eqref{_sum_smallest_Omega_Equation_}
implies that $\Omega^{n-p}$-positivity of 
$R(\eta')$ is equivalent to 
$\alpha_2 + \alpha_4 + ... + \alpha_{2n-2p}\geq 0$.
From \eqref{_in_pairs_Equation_},
this is equivalent to \eqref{_first_2n-2k_eigen_Equation_}.
We proved \ref{_R_eta_Omega_posi_Lemma_}. \endproof

\subsection{$\omega^q$-positive forms in a neighbourhood
  of a subvariety}
\label{_Sibony_Subsection_}

Now we can prove the hypercomplex version of Sibony's
lemma.

\hfill

\theorem\label{_L^1_inte_Theorem_}
Let $M$ be a hyperk\"ahler manifold,
$Z\subset (M, I)$ a compact complex subvariety,
$\codim_\C Z \geq 3$, and 
$\eta \in \Lambda^{2,0}(M\backslash Z,I)$
a real and positive form, which satisfies $\6\eta=0$.
Then $\eta$ is locally integrable everywhere in $M$.

\hfill

{\bf Proof:}
We adapt to hypercomplex situation the coordinate-free
proof of the complex-analytic version of Sibony's lemma,
obtained in \cite{_Verbitsky:omega-psh_}.
In \cite{_Verbitsky:omega-psh_}, the following result was proven.

\hfill

\proposition\label{_forms_zero_in_neighb_Proposition_}
Let $M$ be a K\"ahler manifold, and $Z\subset M$ a complex
subvariety, $\dim_\C Z<p$. Then there exists an open neighbourhood
$U$ of $Z$, and a sequence $\{\rho_i\}$ of 
$\omega^p$-positive, exact, smooth $(1,1)$-forms on $U$ satisfying
the following.
\begin{description}
\item[(i)] For any open subset $V\subset U$,
with the closure $\bar V$ compact and not 
intersecting $Z$, the restriction 
$\rho_i\restrict V$ stabilizes
as $i\arrow \infty$. Moreover, $\rho_i\restrict V$
is strictly $\omega^p$-positive for $i\gg 0$.
\item[(ii)] For all $i$, $\rho_i=0$ in some neighbourhood
  of $Z$.
\item[(iii)] The limit $\rho=\lim\rho_i$ is a strictly
$\omega^p$-positive current on $U$.
\item[(iv)] The forms $\rho_i$ can be written
as $\rho_i = dd^c \phi_i$,
where $\phi_i$ are smooth functions on $U$.
On any compact set not intersecting $Z$,
the sequence $\{\phi_i\}$ stabilizes as  $i\arrow \infty$.
\end{description}

{\bf Proof:} This is \cite{_Verbitsky:omega-psh_}, Proposition 5.3. \endproof

\hfill

We apply \ref{_forms_zero_in_neighb_Proposition_}
to prove \ref{_L^1_inte_Theorem_}.
Let $\phi_i$ be the sequence of functions
defined in a neighbourhood $U \supset Z$
and satisfying conditions of 
\ref{_forms_zero_in_neighb_Proposition_}.
From \ref{_R_differe_Lemma_}, we obtain
\begin{equation} \label{_R_6_bar_6_Equation_}
 R(\6\bar\6\phi_i) = \6 J(\bar\6\phi_i)
\end{equation}
Therefore, $R(\rho_i)$ is $\6$-closed. 
By \ref{_R_eta_Omega_posi_Lemma_},
this form is also $\Omega^{n-1}$-positive.
Since $\eta$ is positive, to show that $\eta$ is
locally integrable on an open set $U\subset M$, 
it suffices to prove that the integral
\begin{equation}\label{_int_w_Omega_univ_bond_Equation_}
\int_D \eta \wedge \Omega^{n-1}\wedge \bar\Omega^n
\end{equation}
is universally bounded, for any compact subset $D\subset U\backslash Z$.
Indeed, 
\[
\int_D \eta \wedge \Omega^{n-1}\wedge \bar\Omega^n 
= \sum_i \int_D \alpha_i \Vol_M
\]
where $\{\alpha_i\}$ are the eigenvalues of $\eta$ considered as 
functions on $M$. 
In \eqref{_int_w_Omega_univ_bond_Equation_},
we may replace $\Omega^{n-1}$ by any 
strictly positive real $(n-1)$-form, and if this integral
us bounded, \eqref{_int_w_Omega_univ_bond_Equation_}
is also bounded. Therefore, \ref{_L^1_inte_Theorem_}
would follow from a universal bound on
\[ \int_D \eta \wedge \rho\wedge \Omega^{n-2}\wedge \bar\Omega^n,\]
where $\rho= \lim R(\rho_i)$ is the form constructed in
\ref{_forms_zero_in_neighb_Proposition_}
(it is smooth outside of $Z$, because $\{\rho_i\}$
stabilizes). Now, a universal bound on
$\int_D \eta \wedge \rho\wedge \Omega^{n-2}\wedge \bar\Omega^n$
would obviously follow from a universal bound on
the integral
\[
\int_D \eta \wedge R(\rho_i)\wedge \Omega^{n-2}\wedge \bar\Omega^n;
\]
this integral is bounded by
\[
\int_U \eta \wedge R(\rho_i)\wedge \Omega^{n-2}\wedge \bar\Omega^n,
\]
because the forms 
$\eta$ and $R(\rho_i)\wedge \Omega^{n-2}$
are positive.\footnote{The product $\eta \wedge R(\rho_i)\wedge
\Omega^{n-2}$  is well defined on the whole
$U$, because $R(\rho_i)$ vanishes in a neighbourhood of $Z$.}

The last integral can be expressed by Stokes' theorem as
\begin{equation}\label{_Stokes_for_neigh_Equation_}
\int_U \eta \wedge R(\rho_i)\wedge \Omega^{n-2}\wedge
\bar\Omega^n = \int_{\6U} 
\eta \wedge J(\bar\phi_i)\wedge \Omega^{n-2}\wedge \bar\Omega^n
\end{equation}
(see \eqref{_R_6_bar_6_Equation_}).
However, the integral $\int_{\6U} 
\eta \wedge J(\bar\phi_i)\wedge \Omega^{n-2}\wedge
\bar\Omega^n$ stabilizes as $i \arrow \infty$,
because $\phi_i$ stabilizes in a neighbourhood of $\6U$.
This shows that \eqref{_int_w_Omega_univ_bond_Equation_}
is universally bounded. We proved 
\ref{_L^1_inte_Theorem_}. \endproof

\section{Skoda-El Mir theorem for hyperk\"ahler manifolds}
\label{_S_E-M_Section_}

We are going to prove  a
hypercomplex analogue of the classical Skoda-El Mir theorem
(\cite{_El_Mir_}, \cite{_Skoda_}, \cite{_Sibony_}, \cite{_Demailly:L^2_}).

\hfill

\definition
Let $M$ be a connected complex manifold, and $Z\subset M$ a closed 
subset. Assume that there exists a nonconstant
plurisubharmonic function $\phi:\; M  \arrow [-\infty, \infty[$,
such that $Z = \phi^{-1}(-\infty)$. Then $Z$ is called
{\bf pluripolar}. 

\hfill

Skoda-El Mir theorem is a result about extending a closed
positive current over a pluripolar set $Z$.

\hfill

\theorem\label{_Sk_E-M_complex_geo_Theorem_}
(\cite{_El_Mir_}, \cite{_Skoda_}, \cite{_Sibony_}, \cite{_Demailly:L^2_})
 Let $X$ be a complex manifold, and $Z$ a closed 
 pluripolar set in $X$. Consider a closed positive current
 $\Theta$ on $X \backslash Z$ which is locally integrable
 around $Z$. Then the trivial extension of $\Theta$ to $X$ is
 closed on $X$.

\endproof

\hfill

The hypercomplex analogue of this theorem goes as follows.

\hfill

\theorem\label{_hk_S-EM_Theorem_}
Let $M$ be a $SL(n, {\Bbb H})$-manifold,
$Z\subset (M,I)$ a 
pluripolar set, and $\eta \in \Lambda^{2p,0}(M\backslash Z, I)$
a form satisfying the following properties.
\begin{description}
\item[(i)] $\eta = J(\bar \eta)$ (reality)
\item[(ii)]  
$\eta(x_1, J(\bar x_1), x_2, J(\bar x_2), ..., x_p, J(\bar x_p))\geq 0$
(weak positivity)
\item[(iii)] $\6\eta=0$ (closedness).
\end{description}
Assume that $\eta$ is integrable around
each point $z\in Z$. Then the trivial extension
of $\eta$ to $M$ is a $\6$-closed $(2p,0)$-current.

\hfill

{\bf Proof:} To prove \ref{_hk_S-EM_Theorem_},
we could repeat the argument proving the Skoda-El Mir theorem
in the hypercomplex setting. However, it is much easier to deduce 
\ref{_hk_S-EM_Theorem_} from the classical Skoda-El Mir. 
Consider the $(p,p)$-form ${\cal R}_{p,p}(\eta)\in \Lambda^{p,p}_I(M)$
obtained as 
\[
{\cal R}_{p,p}(\eta)(x_1, \bar y_1, ..., x_p,  \bar y_p)=
\eta(x_1, J(\bar y_1), ..., x_p, J(\bar y_p)).
\]
where $x_i, y_i\in T^{1,0}(M)$ (see 
Subsection \ref{_p+q,0_and_p,q_explicit_Subsection_}). 

From \ref{_V_main_Proposition_}, it follows that the
$(n+p,n+p)$-form ${\cal R}_{p,p}(\Omega^n) \wedge {\cal R}_{p,p}(\eta)$
is positive in the usual sense if and only if $\eta$ 
is positive in the quaternionic sense,
and closed if and only if $\6\eta=0$.
Now, $\eta$ is closed and positive on $M\backslash Z$, 
hence ${\cal R}_{p,p}(\Omega^n) \wedge {\cal R}_{p,p}(\eta)$ is closed and
positive on $M\backslash Z$ (in the usual sense).
Applying the Skoda-El Mir theorem, we obtain
that a trivial extension of 
${\cal R}_{p,p}(\Omega^n) \wedge {\cal R}_{p,p}(\eta)$ 
is closed on $M$. Applying \ref{_V_main_Proposition_} again,
we find that the trivial extension of $\eta$
to $M$ is $\6$-closed. We proved
\ref{_hk_S-EM_Theorem_}. \endproof

\hfill

{\bf Acknowledgements:}
This article appeared as a byproduct
of a collaboration with Semyon Alesker 
on quaternionic Monge-Ampere equation.

\hfill

\hfill

\small{

\noindent {\sc Misha Verbitsky\\
{\sc  Institute of Theoretical and
Experimental Physics \\
B. Cheremushkinskaya, 25, Moscow, 117259, Russia }\\
\tt verbit@mccme.ru }

}


{\small
\begin{thebibliography}{GMP}

\bibitem[AV]{_Alesker_Verbitsky_HKT_} 
Semyon Alesker, Misha Verbitsky,
{\em Plurisubharmonic functions on hypercomplex manifolds and
 HKT-geometry,} J. Geom. Anal. 16 (2006), no. 3, 375--399.


\bibitem[BS]{_Bando_Siu_} 
 Bando, S., Siu, Y.-T, 
{\it Stable sheaves and Einstein-Hermitian metrics}, 
In: Geometry and Analysis on
Complex Manifolds, Festschrift for Professor S. Kobayashi's 60th Birthday,
ed. T. Mabuchi, J. Noguchi, T. Ochiai, World Scientific, 1994, pp. 39-50.

\bibitem[BS]{_Banos_Swann_}
 Banos, Bertrand; Swann, Andrew;
 {\em Potentials for hyper-K\"ahler metrics with torsion},
arXiv:math/0402366, 
 Classical Quantum Gravity 21 (2004), no. 13, 3127--3135.



\bibitem[BDV]{_BDV:nilmanifolds_}
Maria L. Barberis, Isabel G. Dotti, Misha Verbitsky,
{\em Canonical bundles of complex nilmanifolds, 
with applications to hypercomplex geometry},
arXiv:0712.3863, 22 pages.


\bibitem[Bas]{_Baston_} 
Baston, R. J.
{\it Quaternionic complexes. }
J. Geom. Phys. {\bf 8} (1992), no. 1-4, 29--52. 


\bibitem[Ber]{_Berger:holonomies_}
Berger, M., {\em Sur les groupes d'holonomie des varietes 
a connexion affine et des varietes
riemanniennes}, Bull. Soc. Math. France 83 (1955), 279-330.


\bibitem[Bes]{_Besse:Einst_Manifo_} 
Besse, 
A., {\em Einstein Manifolds}, Springer-Verlag, New York (1987)

\bibitem[B]{_Boyer_} 
Boyer, Charles P.
{\em A note on hyper-Hermitian four-manifolds}.
Proc. Amer. Math. Soc. 102 (1988), no. 1, 157--164. 

\bibitem[Ca]{_Calabi_} 
 Calabi,  E.,
{\em Metriques k\"ahleriennes et fibr\`es holomorphes}, 
Ann. Ecol. Norm. Sup. {\bf 12} (1979), 269-294.  


\bibitem[CS]{_Capria-Salamon_} 
Capria, M. M., Salamon, S. M. 
{\it Yang-Mills fields on quaternionic spaces}, 
Nonlinearity {\bf 1} (1988), no. 4, 517--530. 




\bibitem[D]{_Demailly:L^2_}
Demailly, Jean-Pierre, {\em  $L^2$ vanishing theorems for
  positive line bundles and adjunction theory},
Lecture Notes of a CIME course on "Transcendental Methods
of Algebraic Geometry" (Cetraro, Italy, July 1994),  arXiv:alg-geom/9410022,
and also Lecture Notes in Math., 1646, pp. 1--97, Springer, Berlin, 1996


\bibitem[E]{_El_Mir_}
H. El Mir, {\em Sur le prolongement des courants positifs
fermes,} Acta Math., 153 (1984), 1-45.

\bibitem[GP]{_Gra_Poon_}
Grantcharov, G., Poon, Y. S.,
{\em Geometry of hyper-K\"ahler connections with torsion},
 math.DG/9908015, 
Comm. Math. Phys. 213 (2000), no. 1, 19--37. 



\bibitem[HL1]{_Harvey_Lawson:Psh_}
R. Harvey, B. Lawson, 
{\em  An Introduction to Potential Theory in Calibrated
  Geometry,} arXiv:0710.3920, 45 pages

\bibitem[HL2]{_Harvey_Lawson:Dua_}
R. Harvey, B. Lawson, 
{\em  Duality of Positive Currents and 
Plurisubharmonic Functions in Calibrated Geometry,} arXiv:0710.3921,
29 pages.

\bibitem[HKLR]{_HKLR_}  
N. J. Hitchin, A. Karlhede, 
U. Lindstr\"om, M. Ro\v cek, 
{\em Hyperk\"ahler metrics and supersymmetry}, 
Comm. Math. Phys. {\bf 108} (1987), 535-589.

\bibitem[HP]{_Howe_Papado_}
P.S. Howe, G. Papadopoulos,  {\em Twistor spaces for hyper-K\"ahler
manifolds with torsion} Phys. Lett. B 379 (1996), no. 1-4, 80--86.


\bibitem[H]{_Humphreys_}
J. Humphreys,  {\em Introduction to Lie Algebras and Representation
  Theory},  Graduate Texts
in Mathematics, Springer-Verlag, no. 9, 1972.



\bibitem[LY]{_Leung_} Leung N. C., Yi S.,
{\em Analytic Torsion for Quaternionic manifolds and related topics},
dg-ga/9710022

\bibitem[Ob]{_Obata_} 
Obata, M., {\em Affine connections on manifolds
with almost complex, quaternionic or Hermitian structure}, 
Jap. J. Math., 26 (1955), 43-79.

\bibitem[OSS]{_OSS_} 
 Christian Okonek, Michael 
Schneider, Heinz  Spindler,
{\it Vector bundles on complex projective spaces.}
 Progress in mathematics, vol. 3,
 Birkhauser, 1980.


\bibitem[Sal]{_Salamon_} Salamon, S., {\em Quaternionic Manifolds},
Communicazione inviata all'Instituto nazionale di Alta Matematica
Francesco Severi.


\bibitem[Sib]{_Sibony_}
Sibony, Nessim,
{\em Quelques problemes de prolongement de courants en analyse complexe,}
Duke Math. J. 52, 157-197 (1985).

\bibitem[Sk]{_Skoda_}
H. Skoda,
{\em Prolongement des courants positifs fermes de masse finie,} 
Invent. Math., 66 (1982), 361-376.

\bibitem[V0]{_Verbitsky:Hyperholo_bundles_} 
Verbitsky M., 
{\em Hyperholomorphic bundles over a hyperk\"ahler manifold}, 
alg-geom/9307008, Journ. of Alg. Geom., {\bf 5} no. 4 (1996) pp. 633-669.


\bibitem[V1]{_Verbitsky:Symplectic_II_} 
Verbitsky M., {\em Hyperk\"ahler embeddings and holomorphic 
symplectic geometry II,} alg-geom electronic preprint 9403006 (1994),
14 pages, LaTeX,
also published in: GAFA {\bf 5} no. 1 (1995), 92-104.


\bibitem[V2]{_Verbitsky:hypercomple_}
Verbitsky M.,
{\em Hypercomplex Varieties}, 
alg-geom/9703016,
 Comm. Anal. Geom. {\bf 7} 
(1999), no. 2, 355--396.


\bibitem[V3]{_V:Hyperholo_sheaves_} 
Verbitsky M., {\em
Hyperholomorphic sheaves and new examples of hyperk\"ahler manifolds,}
alg-geom 9712012 - 113 pages, LaTeX 2e; published as a part of a book 
``Hyperk\"ahler manifolds'' (M. Verbitsky, D. Kaledin, 
International Press, Boston).


\bibitem[V4]{_V:reflexive_}
Verbitsky, M.,
{\em Hyperholomorpic connections   on coherent  sheaves 
and stability}, 40 pages, math.AG/0107182

\bibitem[V5]{_Verbitsky:HKT_}
Verbitsky, M., 
{\em 
Hyperk\"ahler manifolds with torsion, supersymmetry and Hodge theory},
math.AG/0112215, 47 pages (Asian J. of Math., 
Vol. 6 (4), December 2002).

\bibitem[V6]{_Verbitsky_HKT-exa_}
M. Verbitsky,
{\em Hyperk\"ahler manifolds with torsion obtained from hyperholomorphic bundles}
 math.DG/0303129, (Math. Res. Lett. 10 (2003), no. 4, 501--513).


\bibitem[V7]{_Verbitsky:canoni_}
M. Verbitsky,  {\em Hypercomplex manifolds with trivial canonical
bundle
 and their holonomy}, 15 pages, arXiv:math/0406537,
``Moscow Seminar on Mathematical Physics, II'',
American Mathematical Society Translations,
{\bf 2}, 221 (2007).


\bibitem[V8]{_Verbitsky:qD_}
M. Verbitsky,
{\em Quaternionic Dolbeault complex and vanishing theorems on
hyperkahler manifolds}, math/0604303, 30 pages, to 
appear in Composition Math.


\bibitem[V9]{_Verbitsky:omega-psh_}
M. Verbitsky,
{\em Plurisubharmonic functions in calibrated geometry and q-convexity},
 arXiv:0712.4036, 27 pages

\bibitem[UY]{_Uhle_Yau_}
 Uhlenbeck K., Yau S. T.,  {\em On the existence of
Hermitian Yang-Mills connections in  stable vector bundles}, 
 Comm. on Pure and Appl. Math., 
{\bf 39}, p. S257-S293 (1986).

\bibitem[Y]{_Yau:Calabi-Yau_} 
Yau, S. T., {\em On the Ricci curvature of a compact K\"ahler manifold 
and the complex Monge-Amp\`ere equation I.}  Comm. on Pure and Appl.
Math. 31, 339-411 (1978).


\end{thebibliography}
}
\end{document}